\input amstex
\documentstyle{amsppt}

\magnification=1200
\parskip 6pt
\pagewidth{5.2in}
\pageheight{7.2in}
\NoRunningHeads
\expandafter\redefine\csname logo\string@\endcsname{}
\NoBlackBoxes

\redefine\frak{\bold}

\redefine\g{\frak g}
\redefine\h{\frak h}
\redefine\p{\frak p}

\redefine\k{\frak k}
\redefine\t{\frak t}
\redefine\n{\frak n}
\redefine\f{\frak f}
\redefine\a{\frak a}
\redefine\u{\frak u}

\define\gc{\frak g^{\C}}
\define\tc{\frak t^{\C}}
\define\kc{\frak k^{\C}}
\define\hc{\frak h^{\C}}

\redefine\Bbb{\bold}
\define\C{\Bbb C}
\define\R{\Bbb R}
\define\Z{\Bbb Z}

\define\N{\Bbb N}
\redefine\H{\Bbb H}

\define\Gc{G^{\C}}

\define\al{\alpha}
\define\be{\beta}
\define\de{\delta}
\define\la{\lambda}

\define\ga{\gamma}

\define\La{\Lambda}
\define\De{\Delta}
\define\Om{\Omega}


\redefine\Im{\operatorname {Im}}

\redefine\dim{\operatorname {dim}}
\define\alg{\operatorname {alg}}
\define\Span{\operatorname {Span}}
\redefine\exp{\operatorname{exp}}
\define\Ad{\operatorname{Ad}}
\define\ad{\operatorname{ad}}

\define\can{\operatorname{can}}
\define\Harm{\operatorname{Harm}}
\define\diag{\operatorname{diag}}

\define\st{\ \vert\ }

\define\ddt{\tfrac d{dt}}

\redefine\sub{\subseteq}

\redefine\ll{\lq\lq}
\redefine\rr{\rq\rq\ }
\define\rrr{\rq\rq}

\redefine\no{\noindent}

\define\sq{{\ssize\sqrt{-1}}\,}
\define\psq{\sq}

\redefine\gal{\g_{\al}}

\topmatter
\title
Harmonic two-spheres in compact symmetric spaces, revisited
\endtitle
\author
F. E. Burstall and M. A. Guest
\endauthor
\endtopmatter

\document

\baselineskip=13pt

\head
Introduction
\endhead

The purpose of this article is to give a new description of
harmonic maps from the two-sphere $S^2$ to a compact 
symmetric space $G/K$, using a method suggested by Morse
theory.  The method leads to surprisingly short proofs
of most of the known results, and also to new results.

The simplest non-trivial examples of such harmonic maps
occur when
$G/K=S^n$ or $\C P^n$. In 1967, E. Calabi gave a construction of
all harmonic maps $S^2\to S^n$. In the early 1980s, all harmonic
maps $S^2\to\C P^n$ were obtained by a similar construction; this
work was carried out independently by various authors, full
details appearing in a 1983 paper of J. Eells and J. C. Wood.
Motivation for continued interest in this problem had 
been provided by mathematical physicists, searching for simple
analogues of the Yang-Mills equations, and indeed the above 
constructions turned out to be related to twistor theory. From 
this point of view, harmonic maps $\phi:S^2\to S^n$ or 
$\C P^n$ were characterized as compositions of the form
$\phi=\pi\circ\Phi$, where $\Phi:S^2\to Z$ is a 
\ll super-horizontal\rr holomorphic map into a \ll twistor
space\rr $Z$, and $\pi:Z\to S^n$ or $\C P^n$ is a
\ll twistor fibration\rrr.  These results are surveyed, with
complete references, in \cite{Ee-Le}.  During the 1980s,
various generalizations of the construction $\phi=\pi\circ\Phi$
appeared, for harmonic maps
$S^2\to G/K$ into various compact symmetric spaces $G/K$.
These were refined into a general and purely Lie-theoretic
treatment in \cite{Bu-Ra}.  We shall refer to this as
\ll the twistor construction\rrr.

Except in the original cases $G/K=S^n$ or $\C P^n$, the
twistor construction does not produce {\it all} harmonic
maps $S^2\to G/K$. Even in the case $G/K=Gr_2(\C^4)$ there
exist harmonic maps which do not arise via the twistor
construction.  The search for a more general construction led
to several special techniques (see the survey \cite{Ee-Le}), by means
of which all harmonic maps $S^2\to Gr_k(\C^n)$ could,
in principle at least, be constructed from 
\ll holomorphic data\rrr. A quite different procedure was
introduced by K. Uhlenbeck in \cite{Uh}, for harmonic maps
$S^2\to U_n$ or $Gr_k(\C^n)$, using ideas from another part
of mathematical physics, namely the theory of integrable
systems. The construction was made more explicit
by Wood in \cite{Wo1} and \cite{Wo2}. However, it remained unclear
how to generalize these results to other compact Lie groups
or symmetric spaces, and even in the case $G=U_n$ or
$G/K=Gr_k(\C^n)$ the results were much less explicit than for
$G/K=S^n$ or $\C P^n$.

In this paper we shall give a straightforward description of
all harmonic maps from $S^2$ to a compact Lie group $G$,
in terms of the special harmonic maps obtained via the
twistor construction. Since any compact symmetric space can be
immersed totally geodesically in its group of isometries,
this description includes all harmonic maps from $S^2$ to 
compact symmetric spaces as well. In the case of an 
{\it inner} symmetric space, we shall describe the harmonic
maps in full detail.

Our description has two main features.
First, we obtain new and very simple proofs of the known
results for $G=U_n$ and $G/K=Gr_k(\C^n)$, and we make those
results more explicit. Second, our method works equally
well for more general $G$ and $G/K$, where few
results were previously known.

The fact that all harmonic maps $S^2\to S^n$ or $\C P^n$
arise from the twistor construction is an easy consequence
of our method. A proof of this fact using
a similar strategy has been given by G. Segal in \cite{Se},
but our approach shows clearly why $S^n$ and $\C P^n$ are
privileged in this way. 
Other well known results on harmonic maps $S^2\to G$
or $G/K$ have simple explanations within our
framework, such as the existence of various 
\ll B\"acklund transformations\rrr, or the factorization 
theorem of \cite{Uh} for
harmonic maps $S^2\to U_n$ or $Gr_k(\C^n)$. We shall give
a proof of this factorization theorem which explains
Lie-theoretically why there is no straightforward
generalization to other Lie groups.

The idea of our method is very simple: given a harmonic map
$\phi$, we apply a (special) family of 
\ll dressing transformations\rrr,
to obtain a family of harmonic maps $\{\phi^t\}_{t\in(0,\infty)}$,
such that $\phi^{\infty}=\lim_{t\to\infty}\phi^t$ is a
harmonic map obtained via the twistor construction.
This leads to explicit results because the family has a
Morse-theoretic interpretation, namely that $\phi^t$ is
obtained by applying the gradient flow of a certain
well known Morse function (to an associated map). 

The main technical result is a classification of harmonic maps into
types, the types (for each fixed $G$ or $G/K$) being indexed
by a finite number of elements in the Lie algebra
$\g$ of $G$. (Strictly speaking, these types correspond
to a finite number of conjugacy classes of parabolic subalgebras
of the complexified Lie algebra $\gc$.) Geometrical 
properties of harmonic maps are reflected by Lie-algebraic 
properties of these elements.  Our first main 
application of this is an
estimate of the \ll minimal uniton number\rr $r(\phi)$ of a 
harmonic map $\phi$.  For $G=U_n$, Uhlenbeck showed that
the maximal value of $r(\phi)$ is $n-1$. We obtain the
maximal values for all $G$ and $G/K$. The second main application
is a \ll Weierstrass formula\rrr, by means of which 
harmonic maps $S^2\to G$ or $G/K$ are described explicitly
in terms of meromorphic functions on $S^2$ (i.e. rational
functions). This includes the explicit formulae which were
known for $G/K=S^n$ or $\C P^n$, and underlies 
Wood's formulae for $G=U_n$.  
It generalizes a Weierstrass formulae of R. Bryant for harmonic maps 
obtained via the twistor construction (see \cite{Br1},\cite{Br2}).

While this method is quite succesful in extending and
clarifying the existing theory, we believe that our
point of view can be justified further.  For example,
our Weierstrass formulae are sufficiently more explicit
than the existing formulae (even in the case $G=U_n$ or
$G/K=Gr_k(\C^n)$) that they can be used to study the
space \ll Harm\rr of harmonic maps. Similar Morse-theoretic ideas
have already been used to study $\Harm(S^2,S^n)$ and
$\Harm(S^2,\C P^n)$ in \cite{Gu-Oh} and \cite{Fu-Gu-Ko-Oh},
and our methods should permit generalizations to other
$G$ and $G/K$. Secondly, our approach to Weierstrass formulae
was influenced by the paper \cite{Do-Pe-Wu}, in which
a general scheme is suggested for constructing harmonic
maps $M\to G$ or $G/K$ in terms of \ll Weierstrass data\rrr,
for compact Riemann surfaces $M$ of arbitrary genus $g$.  Our
results show how this scheme may be implemented in the case 
$g=0$. The case $g=1$ has received much recent attention,
and the scheme seems likely to be implementable also in
this case (see \cite{Bu-Fe-Pe-Pi}, \cite{Bu}, \cite{Mc}).

The paper is arranged as follows.
In \S 1 we review the standard reformulation of harmonic maps
in terms of extended solutions, from \cite{Uh}. We also
discuss complex extended solutions, a notion suggested in
\cite{Do-Pe-Wu}. We give a general definition of uniton number
and minimal uniton number.  Morse theory enters the picture
in the form of the classical energy functional on the loop
group $\Om G$. This is well known from \cite{Bo}, and goes 
back to Morse himself. We use the modern version provided by
the theory of loop groups (see \cite{Pr-Se}), and so we review this
in \S 2.   The main ingredients are the Birkhoff and Bruhat 
decompositions, which are most conveniently described in a 
purely algebraic manner. In \S 3 we review the twistor construction
for harmonic maps into symmetric spaces, following \cite{Bu-Ra}.
Such harmonic maps are characterized by the fact that they
correspond to
$S^1$-invariant extended solutions, for a certain $S^1$-action
on $\Om G$. (This relevance of this action for harmonic maps 
was first pointed out by C.-L. Terng.) These
harmonic maps are fundamental in our theory.  As explained earlier,
our method is to analyse general extended solutions in terms of  
associated $S^1$-invariant extended solutions.  This is where
Morse theory, in the guise of the Bruhat decomposition, enters.
Our results are given in \S 4 for harmonic maps $S^2\to G$,
and in \S 5 for harmonic maps $S^2\to G/K$.
In fact, the results apply equally well to harmonic maps 
$M\to G$ or $G/K$ of finite
uniton number, where $M$ is any Riemann surface.
This includes those maps described in the literature as
isotropic, or superminimal. Finally, there are two short
appendices.  In Appendix A we collect together the special
results which apply to harmonic maps of \ll low uniton number\rrr.
To some extent this explains the historical development
of the subject, in which the fundamental results for
$S^n$ and $\C P^n$ were gradually extended to more complicated
symmetric spaces.  Appendix B
considers the role of the Birkhoff decomposition, and the
relation between our method and the \ll meromorphic potentials\rr
of \cite{Do-Pe-Wu}.

The second author was partially supported by the U.S. 
National Science Foundation.

\newpage
\head
\S 1  Review of extended solutions
\endhead

\subheading{Extended solutions and harmonic maps}

Let $G$ be a connected compact Lie group, with Lie algebra $\g$.
Let $M$ be a connected (but not necessarily compact)
Riemann surface, without boundary.  To study harmonic
maps $M\to G$, we use the well known correspondence between
harmonic maps and extended solutions, as in \cite{Uh}.
The relevant definitions will be summarized briefly
in this section.

If $X$ is any manifold, we denote the 
\ll free\rr and \ll based\rr
loop spaces of $X$ (respectively) by $\La X$, $\Om X$,
where
$$\align
\La X&=\{\ga:S^1\to X \st \ga \ \text{is smooth}\}\\
\Om X&=\{\ga\in\La X\st \ga(1)=x\}
\endalign
$$
where $x$ is a fixed basepoint of $X$. If $X$ is a
group, we refer to these as the free and based 
{\it loop groups} of $X$, and we take $x=e$, the identity 
element of the group.
We consider $S^1$ here to be the set of unit complex
numbers, i.e. $S^1=\{\la\in\C\st \vert\la\vert=1\}$.

Let $\Gc$ be the complexification of $G$, with Lie
algebra $\gc$ (thus $\gc=\g\otimes\C$).

\no{\it Definition:}
A map $\Phi:M\to\Omega G$ is an {\it extended solution}
if it satisfies the equation
$$
\Phi(z,\la)^{-1}\Phi_{z}(z,\la)=
\tfrac12(1-\tfrac1\la)A(z)
\tag{$E$}
$$
for some map $A:M\to\gc$.
(Here we write $\Phi(z,\la)$ for $\Phi(z)(\la)$, 
with $z\in M$ and $\la\in S^1$, and we write $\Phi_z$
for the partial derivative of $\Phi$ with respect to $z$.)

\proclaim{Theorem 1.1  \cite{\bf Uh} } 

\no(1) If $\Phi:M\to\Om G$ is an extended solution, then the
map $\phi:M\to G$ defined by $\phi(z)=\Phi(z,-1)$
is harmonic.

\no(2) If $\phi:M\to G$ is harmonic, and if $M$ is
simply connected, then there exists an extended
solution $\Phi:M\to\Om G$ such that $\phi(z)=\Phi(z,-1)$.
This $\Phi$ is unique up to multiplication on the left
by an element $\ga\in\Om G$ such that $\ga(-1)=e$.
\qed
\endproclaim

It is sometimes convenient to say that an extended solution
$\tilde\Phi:M\to\La G$ is a map which satisfies the following
equation:
$$
\tilde\Phi(z,\la)^{-1}\tilde\Phi_{z}(z,\la)=
B(z)+\tfrac1\la C(z)
\tag{$\tilde E$}
$$
for some maps $B,C:M\to\gc$. Certainly a solution of
$(E)$ is a solution of $(\tilde E)$; conversely, if
$\tilde\Phi$ is a solution of $(\tilde E)$, then
it is easy to verify
that $\Phi(z,\la)=\tilde\Phi(z,\la)\tilde\Phi(z,1)^{-1}$
is a solution of $(E)$.

In this paper we shall be concerned entirely with
extended solutions of \ll finite uniton number\rrr. 
This concept was introduced in \cite{Uh} for the
case $G=U_n$.  We shall give a definition for general
$G$, which extends the definition of \cite{Uh} in a
natural way.  First, recall that the adjoint representation
is the homomorphism
$$
\Ad:G\to O(\g),\quad \Ad(g)\xi = 
\ddt g (\exp\,t\xi) g^{-1}\vert_0,
$$
where $O(\g)$ is the orthogonal group of (the vector
space) $\g$ with respect to a fixed bi-invariant inner
product on $\g$. If $G$ is a matrix group, then 
$\Ad(g)\xi=g\xi g^{-1}$. The image $\Ad G$ is isomorphic to $G/Z(G)$,
where $Z(G)$ is the centre of $G$.

\no{\it Definition:}
A loop $\ga:S^1\to G$ is {\it algebraic} if
it is the restriction of an algebraic map
$\C^\ast\to\Gc$.

\no The algebraic loops in $\Om G$ form a subgroup, 
$\Om_{\alg}G$. If $\theta:G\to U_n$ is a representation
of $G$, and $\ga\in\Om_{\alg}G$, then $\theta(\ga)$ is
necessarily a polynomial in $\la$ and $\la^{-1}$.
Taking $\theta=\Ad$, we obtain a filtration
$$
\{e\}=\Om_{\alg}^0 G\sub \Om_{\alg}^1 G\sub \dots 
\sub\cup_{k\ge 0}\Om_{\alg}^k G=\Om_{\alg}G
$$
by defining $\Om_{\alg}^k G$ to be the set of loops $\ga$
such that $\Ad(\ga)$ is of the form 
$\sum_{\vert i\vert\le k} \la^i T_i$.

\no{\it Definition:}
Let $\Phi:M\to \Om G$ be an extended solution.
If  $\Phi(M)\sub \Om_{\alg}^k G$ and
$\Phi(M)\not\sub \Om_{\alg}^{k-1} G$, 
we say that $\Phi$ has {\it uniton number} $k$. We
write $r(\Phi)=k$.

That this is not a vacuous definition is shown by
the next theorem:

\proclaim{Theorem 1.2  \cite{\bf Uh},\cite{\bf Se} }
Let $\Phi:M\to \Om G$ be an extended solution. If
$M$ is compact, then there exists some $\ga\in \Om G$
and some $k\ge 0$ such that $\ga\Phi(M)\sub \Om_{\alg}^k G$.
\qed
\endproclaim

\no Actually, Uhlenbeck only considers the case $G=U_n$, and proves
that if $\Phi:M\to\Omega U_n$ is an extended solution
then there is some $\ga\in\Om U_n$ such that $\ga\Phi$ has
the form $\sum_{i=0}^k \la^i A_i$. It is easy to deduce
Theorem 1.2 from this, by viewing $\Ad G$ as a subgroup
of the unitary group $U(\gc)$.
Uhlenbeck goes on to define the uniton number of 
$\ga\Phi=\sum_{i=0}^k \la^i A_i$ to be $k$. This definition
coincides with ours so long as $A_0\ne 0$. (Our definition 
factors out the effect of scalar loops $\la^k$, which 
are assigned uniton number $0$.)

When $M$ is both compact and simply connected, i.e. $M=S^2$,  
Theorems 1.1 and 1.2 say that any harmonic map $M\to G$
corresponds to an extended solution of finite uniton number.
Because of this, our results will apply to arbitrary harmonic
maps $S^2\to G$, but only to those harmonic maps $M\to G$
which come from extended solutions of finite uniton number. 

Since the correspondence between harmonic maps and extended
solutions is not one to one (even when $M$ is simply
connected), we need a separate definition for the \ll uniton
number of a harmonic map\rrr. Let $\phi:M\to G$ be a harmonic map,
and let $\Phi:M\to \Om G$ be an extended solution with 
$r(\Phi)<\infty$ such that $\Phi(z,-1)=\phi(z)$. Associated
to $\phi$ we have the non-negative integer
$$
\tilde r(\phi)=\min\{ r(\ga\Phi)\st \ga\in \Om_{\alg}G \}.
$$
Unlike the definition of $r(\Phi)$, however, the definition of 
$\tilde r(\phi)$ depends on the local isomorphism class of $G$.
This is not appropriate for our purposes, as harmonic maps
into locally isomorphic Lie groups are (locally) equivalent.
So we modify $\tilde r(\phi)$ as follows:

\no{\it Definition:}
Let $\phi:M\to G$ be a harmonic map.
Assume that there exists an extended solution
$\Phi:M\to \Om G$ such that $r(\Phi)<\infty$ 
and $\Phi(z,-1)=\phi(z)$. We define the
{\it minimal uniton number} of $\phi$ as the non-negative
integer
$$
r(\phi)=\min\{ r(\ga \Ad\Phi) \st \ga \in \Om_{\alg}\Ad G\}.
$$

\no The minimal uniton number of $\phi$, which depends 
only on $\Ad\phi$, is a
measure of the \ll complexity\rr of the harmonic map.
In the case $G=U_n$, our
definition agrees exactly with Uhlenbeck's definition.

We conclude with a simple example which illustrates the
advantages of $r(\phi)$ over $\tilde r(\phi)$. Let
$\phi:S^2\to S^2$ be a holomorphic map. Consider the
totally geodesic embedding $i:S^2\cong \C P^1\to SU_2$, 
$l\mapsto (p-p^\perp)(\pi_l-\pi_l^\perp)$, where $\pi_l$
denotes Hermitian projection on $l$, and $p$ is a fixed
Hermitian projection operator of rank $1$. The composition
$i\circ\phi$ is a harmonic map $S^2\to SU_2$. A suitable
extended solution $\Phi:S^2\to\Om SU_2$ is given by
$\Phi(z,\la)=(p+\tfrac1\la p^\perp)
(\pi_{\phi(z)}+\la \pi_{\phi(z)}^\perp)$. It is easy to see
that $\tilde r(\phi)=r(\Phi)=2$. However, we find that $r(\phi)=1$, 
because $\Om PSU_2=\Om SU_2\sqcup [p+\tfrac1\la p^\perp]
\Om SU_2$. (We consider $\Om SU_2$ as the identity component
of $\Om PSU_2$, and $[p+\tfrac1\la p^\perp]$ as the
topologically non-trivial loop in $PSU_2\cong U_2/Z(U_2)$
which is given by the loop $p+\tfrac1\la p^\perp$ in
$U_2$.)  Thus, it is only by passing from $G$ to $\Ad G$
that we are able to remove the superfluous factor 
$p+\tfrac1\la p^\perp$.

\subheading{Complex extended solutions}

It is well known (see \cite{Pr-Se}) that the loop group
$\Om G$ admits the structure of an (infinite dimensional)
complex manifold. One way to describe this complex structure
comes from the \ll Iwasawa decomposition\rr of
the loop group $\La \Gc$.  In the following statement,
$\La^+\Gc$ denotes the subgroup of $\La \Gc$ consisting
of maps $S^1\to\Gc$ which extend holomorphically to the region
$D^+=\{\la\st\vert\la\vert<1\}$: 

\proclaim{Theorem 1.3 \cite{\bf Pr-Se} }
The product map
$\La^+\Gc\times\Om G\to\La\Gc$ is a diffeomorphism.
In particular, $\La\Gc=\Om G\,\La^+\Gc$, and any
$\ga\in\La\Gc$ may be written uniquely in the form
$\ga=\ga_u\ga_+$, where $\ga_u\in\Om G$
and $\ga_+\in \La^+\Gc$.
\qed
\endproclaim

This result has a number of useful consequences.
The first of these concerns the definition of the
complex structure on $\Om G$.
Each of the groups $\La\Gc$, $\La^+\Gc$ is, in a
natural way, a complex Lie group, and so the homogeneous space 
$\La\Gc/\La^+\Gc$ is a complex manifold.
By the theorem, we
have an identification $\Om G\cong \La\Gc/\La^+\Gc$, so
this gives the required complex structure on 
$\Om G$. It follows from standard theory of homogeneous spaces
that the holomorphic tangent bundle 
$T^{1,0}\Om G$ of $\Om G$ may be
identified with $\La\Gc\times_{\La^+\Gc}(\La\gc/\La^+\gc)$.

Second, the theorem shows that there
is a natural action of the complex group $\La\Gc$ on
the complex manifold $\Om G$, given by the natural
action of $\La\Gc$ on $\La\Gc/\La^+\Gc$. We can
express this action in the following way: if 
$\ga\in\La\Gc$ and $\de\in\Om G$, then
$$
\ga\cdot\de = (\ga\de)_u.
$$
It turns out that $\La\Gc$ acts as a symmetry group
of the extended solution equation $(E)$, i.e. if $\Phi$ is 
an extended solution then so is $\ga\cdot\Phi$. It was
shown in \cite{Gu-Oh} that this action is essentially
the same as the \ll dressing action\rr introduced in \cite{Uh}.

The third consequence is that we may reformulate the 
extended solution equation (E) in terms of the complex loop group.
Let us write $\Phi=[\Psi]$, where $\Psi:M\to\La\Gc$.
(Such a map $\Psi$ exists locally, at least. If $M$ is simply
connected, $\Psi$ exists globally.) Note that $[\Psi]$
may also be written $\Psi_u$.

\proclaim{Proposition 1.4}
The map $\Phi$ satisfies $(E)$ if and only if $\Psi$
satisfies
$$
\cases
\Im\,\la \Psi^{-1}\Psi_z&\sub \ \La^+\gc \\
\Im\,\Psi^{-1}\Psi_{\bar z}&\sub \ \La^+\gc
\endcases
\tag{$E^{\C}$}
$$
(where \ll\ $\Im$\rr denotes \ll image\rrr).
\endproclaim

\demo{Proof} By Theorem 1.3
we may write $\Psi=\Phi\De$, where $\De:M\to\La^+\Gc$.
The identities
$\Psi^{-1}\Psi_z=\De^{-1}\Phi^{-1}\Phi_z\De + \De^{-1}\De_z$,
$\Psi^{-1}\Psi_{\bar z}=
\De^{-1}\Phi^{-1}\Phi_{\bar z}\De + \De^{-1}\De_{\bar z}$
show that $(E^{\C})$ is equivalent to
$$
\cases
\Im\,\la \Phi^{-1}\Phi_z&\sub \ \La^+\gc \\
\Im\,\Phi^{-1}\Phi_{\bar z}&\sub \ \La^+\gc
\endcases.
\tag{$E'$}
$$
But $(E')$ is equivalent to $(E)$, because
if $\Phi^{-1}\Phi_z=\sum \la^i A_i$ then
$\Phi^{-1}\Phi_{\bar z}=\sum \la^{-i} \bar A_i$,
where the bar denotes complex conjugation in $\gc$
with respect to the real form $\g$.
\qed\enddemo

Because of the identification 
$T^{1,0}\Om G\cong \La\Gc\times_{\La^+\Gc}(\La\gc/\La^+\gc)$,
the condition $\Im\,\Psi^{-1}\Psi_{\bar z}$ {} $\sub \La^+\gc$
says that $\Phi=[\Psi]:M\to \Om G$ is holomorphic.
This is the key to yet another reformulation of the
extended solution equation. Let 
$\La^{\ast}\Gc$ denotes the subgroup of $\La \Gc$ consisting
of maps $S^1\to\Gc$ which extend holomorphically to the region
$D^{\ast}=\{\la\st 0<\vert\la\vert<1\}$.

\no{\it Definition:}
A holomorphic map $\Psi:M\to\La^{\ast}\Gc$ is a 
{\it complex extended solution}
if $\Im\,\la\Psi^{-1}\Psi_z\sub\La^+\gc$.

\proclaim{Theorem 1.5  \cite{\bf Do-Pe-Wu} }

\no(1) If $\Psi:M\to\La^{\ast}\Gc$ is a complex extended solution, then
$\Phi=\Psi_u$ is an extended solution.

\no(2) If $\Phi:M\to\Om G$ is an extended solution,
and $z_0$ is any point of $M$, 
then there exists a neighbourhood $M_0$ of $z_0$
and a complex extended solution
$\Psi:M_0\to\La\Gc$ such that $\Phi\vert_{M_0}=\Psi_u$.
\qed
\endproclaim

\no We shall
give a proof of Theorem 1.5 in the case of 
extended solutions of finite uniton number, in \S 4.

\newpage
\head
\S 2 Review of loop groups
\endhead

\subheading{Birkhoff and Bruhat decompositions}

The (based) loop group $\Om G$ may be identified with the
complex homogeneous space $\La\Gc/\La^+\Gc$, as we remarked in
\S 1. This is proved in \cite{Pr-Se} by showing that
both spaces may be identified with a certain infinite
dimensional complex Grassmannian.  The
analogy between $\Om G$ and (finite dimensional)
Grassmannians $Gr_k(C^n)$ is in fact one of the main themes
of \cite{Pr-Se}. In this section we shall review one
aspect of this, namely the Birkhoff and Bruhat decompositions
of $\Om G$. These are analogous to the Schubert decomposition
of $Gr_k(C^n)$.

Let us choose a maximal torus $T$ of $G$. The group of
homomorphisms $S^1\to T$ may be identified 
with the integer lattice $I=(\exp 2\pi)^{-1}(e)\cap\t$
in $\t$, by associating to $\xi\in I$ the homomorphism
$\ga_{\xi}:\la=e^{\psq t}\mapsto \exp\,t\xi$.
Let us also choose a fundamental Weyl chamber in $\t$.
The intersection of $I$ with this will be denoted $I'$.
Whereas $I$ parametrizes homomorphisms $S^1\to T$, $I'$
parametrizes conjugacy classes of homomorphisms $S^1\to G$
(or, equivalently, orbits of homomorphisms $S^1\to T$
under the Weyl group).

\proclaim{Theorem 2.1 {\bf (Chapter 8 of \cite{\bf Pr-Se})}  }

\no (1) Birkhoff decomposition:
$
\La\Gc=\bigsqcup_{\xi\in I'}\La^-\Gc\,\ga_{\xi}\,\La^+\Gc.
$

\no (2) Bruhat decomposition:
$
\La_{\alg}\Gc=\bigsqcup_{\xi\in I'}
\La^+_{\alg}\Gc\,\ga_{\xi}\,\La^+_{\alg}\Gc.
$
\endproclaim

\no (A loop $\ga\in\La\Gc$ is said to be algebraic if both
$\Ad\ga$ and $\Ad\ga^{-1}$ are of the form 
$\sum_{\vert i\vert\le k} \la^i T_i$; $\La_{\alg}\Gc$ denotes
the subgroup of algebraic loops, and $\La^+_{\alg}\Gc$ denotes
$\La_{\alg}\Gc\cap\La^+\Gc$.)

\proclaim{Corollary 2.2 {\bf (Chapter 8 of \cite{\bf Pr-Se})} }

\no (1) Birkhoff decomposition:
$
\Om G=\bigsqcup_{\xi\in I'}\La^-\Gc\cdot \ga_{\xi}.
$

\no (2) Bruhat decomposition:
$
\Om_{\alg} G=\bigsqcup_{\xi\in I'}
\La^+_{\alg}\Gc\cdot\ga_{\xi}.
$
\endproclaim

\no The nature of these decompositions is best understood in terms
of Morse theory.  We shall review this next, 
following \cite{Pr}, \cite{Pr-Se}.

\subheading{Morse theory}  

Let $E:\Om G\to\R$ denote the usual energy functional
on paths, $E(\ga)=\int_{S^1}\vert\ga'\vert^2$. 
This is a Morse-Bott function, i.e. each of its critical
manifolds is non-degenerate. The critical points are
the geodesics in $G$ which pass through the identity element,
i.e. the homomorphisms $S^1\to G$, and the (connected) 
critical manifolds are the conjugacy classes of such 
homomorphisms. We write
$$
\Om_{\xi}=\{g\ga_{\xi} g^{-1}\st g\in G\}
$$
for the conjugacy class of the homomorphism $\ga_{\xi}$.
We write
$$
\align
S_{\xi}&=\{\ga\in\Om G\st 
\ \text{$\ga$ flows into $\Om_{\xi}$}\}\\
U_{\xi}&=\{\ga\in\Om G\st 
\ \text{$\ga$ flows out of $\Om_{\xi}$}\}
\endalign
$$
for the unstable and stable manifolds of $\Om_{\xi}$
(where \ll flow\rr refers to the flow of the gradient
vector field $-\nabla E$). The unstable manifold
$U_{\xi}$ has the structure of a vector bundle 
over $\Om_{\xi}$; we denote the bundle map by 
$u_{\xi}:U_{\xi}\to\Om_{\xi}$. The rank of this bundle is the
Morse index of $\ga_{\xi}$, which may be expressed
in terms of the roots of $G$ (we shall give the
formula shortly).  
All these facts were discovered by Bott (\cite{Bo}).

The relation with the Birkhoff and Bruhat decompositions
is given by the following result of Pressley:

\proclaim{Theorem 2.3  \cite{\bf Pr}}

\no (1) $S_{\xi}=\La^-\Gc\cdot \ga_{\xi}$.

\no (2) $U_{\xi}=
\La^+_{\alg}\Gc\cdot\ga_{\xi}$.
\endproclaim

\no (A similar theorem holds for the Schubert decomposition 
of the Grassmannian $Gr_k(\C^n)$, and, in fact, for any 
generalized flag manifold - see \cite{Pa}.) 

The theorem is proved in \cite{Pr} by taking a
\ll Hamiltonian\rr point of view, in contrast to the
differential geometric methods of \cite{Bo}. Namely, the
energy functional $E$ is viewed as a Hamiltonian function
associated to a symplectic action of the group $S^1$ on
$\Om G$; this action is given by
$$
u\cdot\ga(\la)=\ga(u\la)\ga(u)^{-1},
\quad u\in S^1,\ga\in\Om G.
$$
The critical points of $E$ are the fixed points of this
action.  
The symplectic structure of $\Om G$ actually
comes from a K\"ahler structure.  If the gradient
is taken with respect to the K\"ahler metric, then the
flow of $-\nabla E$ can be described {\it explicitly} 
in terms of
the \ll complexification\rr of the $S^1$-action.
Let $\C^\ast_{{\sssize\ge} 1}=\{\la\in\C\st 1\le\vert\la\vert<\infty\}$.
This semigroup acts on $\Om G$ by
$$
u\cdot\ga(\la)=\ga(u\la)\La^+\Gc\in \La\Gc/\La^+\Gc,
\quad u\in \C^\ast_{{\sssize\ge} 1},\ga\in\Om G.
$$
It turns out that the flow line of 
$-\nabla E$ starting
at a point $\ga\in\Om G$ is given by the action of
the subsemigroup $[1,\infty)$ of $\C^\ast_{{\sssize\ge} 1}$.
The flow line of $\nabla E$ starting at $\ga$
is defined for all time if and only $\ga\in\Om_{\alg}G$,
in which case it is given by the action of the
subsemigroup $(0,1]$ of 
$\C^\ast=\{\la\in\C\st 0<\vert\la\vert<\infty\}$.
Note that the formula for the action of $\C^\ast_{{\sssize\ge} 1}$
on $\Om G$ extends to an action of $\C^\ast$ on $\Om_{\alg}G$.
It was observed by Terng (see \S 7 of \cite{Uh}) that this
action preserves the extended solution equation $(E)$.
(In combination with the earlier action of $\La_{\alg}\Gc$,
this gives a symmetry group 
$\C^\ast\ltimes\La_{\alg}\Gc$ of $(E)$).

\subheading{The holomorphic vector bundle 
$u_{\xi}:U_{\xi}\to\Om_{\xi}$}

We shall need an explicit description of the vector bundle
$u_{\xi}:U_{\xi}\to\Om_{\xi}$.
Let $\De$ be the set of roots of $\gc$ with respect to the
maximal torus $T$; thus $\De\sub\sq\t^\ast$. We have already
chosen a fundamental Weyl chamber, so we have a decomposition
$\De=\De^+\sqcup\De^-$ of $\De$ into the subsets of positive
and negative roots. As
above, we take $\xi\in I'$. It follows that 
$\al(\xi)/\psq\in\Z$ for any $\al\in\De$, and
$\al(\xi)/\psq\ge 0$ for any $\al\in\De^+$.
Let
$\gal$ be the root space of $\al$ (thus, 
$\ad(\tau)\eta=\al(\tau)\eta$
for all $\tau\in\t$, $\eta\in\gal$). Let $\g^{\xi}_i$ be the
$\psq i$-eigenspace of $\ad\xi$. With these definitions, we have:
$$
\gc={\tsize\bigoplus}_{i}\g^{\xi}_i,\quad\quad
\g^{\xi}_i={\tsize\bigoplus}_{\al(\xi)=\psq i}\gal.
$$

\no{\it Definition:} The {\it height} of $\xi$ is the
non-negative integer 
$r(\xi)=\max\{i\st\g^{\xi}_i\ne 0\}$.

First we shall describe $\Om_{\xi}$ abstractly, as a
complex homogeneous space. Consider the orbit $\Gc\cdot\ga_{\xi}$
in $\La\Gc\cdot\ga_{\xi}=\Om G$.
The isotropy subgroup at $\ga_{\xi}$ is the subgroup
$P_{\xi}=\Gc\cap \ga_{\xi} (\La^+\Gc) \ga_{\xi}^{-1}$
of $\Gc$.  The Lie algebra of $P_{\xi}$ is
$$
\align
\p_{\xi}&=\gc\cap \Ad(\ga_{\xi})\La^+\gc\\
&=\tc\oplus
(
{\tsize\bigoplus}_{\al\in\De^-}\gal
)
\oplus
(
{\tsize\bigoplus}_{\al\in\De^+,\al(\xi)=0}\gal
)
\\
&={\tsize\bigoplus}_{i\le 0}\g^{\xi}_i.
\endalign
$$
Here we use the fact that 
$\Ad\,\ga_{\xi}=\Ad\,\exp\,t\xi=e^{\ad t\xi}$, which
is given by multiplication by $\la^i$ on $\g^{\xi}_i$.
Evidently we have $G\cdot\ga_{\xi}\sub \Gc\cdot\ga_{\xi}$.
It turns out, however, that these are equal.
This follows from the
Iwasawa decomposition $\Gc=G A N$, where $A,N$ are the
connected subgroups of $\Gc$ corresponding to the Lie
algebras 
$\a=\sq\t$, $\n={\tsize\bigoplus}_{\al\in\De^-}\,\gal$.
(Note that $AN\sub P_{\xi}$.)
We therefore have an explicit identification
$$
\Om_{\xi}\cong \Gc/P_{\xi},
$$
and the holomorphic tangent bundle of $\Om_{\xi}$
is
$T^{1,0}\Om_{\xi}\cong \Gc\times_{P_{\xi}}(\gc/\p_{\xi})$.

Next, we shall describe $U_{\xi}$ in a similar fashion.  
It is by definition a complex homogeneous space of the group 
$\La_{\alg}^+\Gc$, and the isotropy subgroup at $\ga_{\xi}$ 
is the subgroup
$
\La_{\alg}^+\Gc\cap\ga_{\xi}(\La^+\Gc)\ga_{\xi}^{-1}.
$
The Lie algebra of the isotropy subgroup is
$$
\align
\La_{\alg}^+\gc\cap \Ad(\ga_{\xi})\La^+\gc
&= {\tsize\bigoplus}_{i\ge0,\,\al(\xi)/\psq\le i}\la^i\gal\\
&= {\tsize\bigoplus}_{i\ge0}\la^i
(
{\tsize\bigoplus}_{j\le i}\g^{\xi}_j
)
\\
&= {\tsize\bigoplus}_{0\le i\le r(\xi)}\la^i\f^{\xi}_i
\endalign
$$
where $\f^{\xi}_i=\bigoplus_{j\le i}\ \g^{\xi}_j$.

Finally, we claim that the bundle $u_{\xi}:U_{\xi}\to\Om_{\xi}$
is given simply by the natural map
$$
\La_{\alg}^+\Gc/
\La_{\alg}^+\Gc\cap \ga_{\xi} (\La^+\Gc) \ga_{\xi}^{-1}
\to
\Gc/P_{\xi}
$$
which is induced by the homomorphism
$\La_{\alg}^+\Gc\to \Gc,\quad \ga\mapsto \ga(0)$.
This follows from the fact that the flow of $\nabla E$ is
given by the action of $(0,1]$, as described earlier.

We note that this is in fact a {\it holomorphic}
fibre bundle. To see that it is a {\it vector} bundle (and
to calculate its rank), we consider the fibre over 
$\ga_{\xi}$. From the above description, this fibre is
$\La_{e,\alg}^+\Gc\cdot\ga_{\xi}$, where
$\La_{e,\alg}^+\Gc=\{\ga\in\La_{\alg}^+\Gc\st \ga(0)=e\}$.
It may therefore be identified with the homogeneous space
$$
\La_{e,\alg}^+\Gc\ /\ 
\La_{e,\alg}^+\Gc\cap \ga_{\xi} (\La^+\Gc) \ga_{\xi}^{-1}.
$$ 
This is a complex manifold of dimension
$\sum_{\al\in\De^+,\al(\xi)\ne 0}(\al(\xi)/\psq-1)$
(i.e. the dimension of 
$\bigoplus_{0<i<\al(\xi)/\psq}\la^i\gal$).
We can be even more explicit:

\proclaim{Proposition 2.4} The fibre of 
$u_{\xi}:U_{\xi}\to\Om_{\xi}$ over $\ga_{\xi}$ is 
$$
\La_{e,\alg}^+\Gc\cdot\ga_{\xi}=
\exp\,\u_{\xi}\cdot \ga_{\xi} \cong \exp\,\u_{\xi} \cong \u_{\xi}
$$
where $\u_{\xi}$ is the (finite dimensional) nilpotent
subalgebra of $\La_{\alg}^+\gc$ defined by
$$
\u_{\xi} = {\tsize\bigoplus}_{0<i<r(\xi)}\la^i(\f^{\xi}_i)^\perp,
\quad\quad (\f^{\xi}_i)^\perp={\tsize\bigoplus}_{i<j\le r(\xi)}\g^{\xi}_j.
$$
\endproclaim

\demo{Proof} By construction, the isotropy subgroup of
$\exp\,\u_{\xi}$ at $\ga_{\xi}$ is trivial, so we have
$\exp\,\u_{\xi}\cdot \ga_{\xi} \cong \exp\,\u_{\xi}$. 
Since $\u_{\xi}$ is a finite dimensional nilpotent
Lie algebra, its exponential map is biholomorphic,
so we have $\exp\,\u_{\xi} \cong \u_{\xi}$.
It remains to show that 
$\exp\,\u_{\xi}\cdot \ga_{\xi}=\La_{e,\alg}^+\Gc\cdot\ga_{\xi}$.
Since $\exp\,\u_{\xi}$ is a subgroup of $\La_{e,\alg}^+\Gc$,
it suffices to show that both orbits have the
same dimension (cf. \cite{Pr}, page 558, conditions (A) and (B)).
That they do follows from the formula above for the
dimension, and the definition of $\u_{\xi}$.
\qed\enddemo

\no We deduce from this that $u_{\xi}:U_{\xi}\to\Om_{\xi}$
is a holomorphic vector bundle, of rank 
\newline
$\sum_{\al\in\De^+,\al(\xi)\ne 0}(\al(\xi)/\psq-1)$.

For later use, we shall compute explicitly the inclusions
of holomorphic tangent spaces
$$
\CD
T^{1,0}\Om_{\xi} @>>> T^{1,0}U_{\xi} @>>> T^{1,0}\Om G.
\endCD
$$
It suffices to work over the point $\ga_{\xi}$. Applying
left translation by $\ga_{\xi}^{-1}$, we have a sequence
$$
\CD
\gc/\p_{\xi} @>{i_1}>> \La_{\alg}^+\gc/
\La_{\alg}^+\gc\cap \Ad(\ga_{\xi})\La^+\gc
@>{i_2}>> \La\gc/\La^+\gc.
\endCD
$$

\proclaim{Lemma 2.5}
The maps $i_1$, $i_2$ are given by
$$
\align
i_2\circ i_1[\eta]&=[\la^{-i}\eta]\\
i_2[\la^j\eta]&=[\la^{j-i}\eta]
\endalign
$$
where $\eta\in\g^{\xi}_i$.
\endproclaim

\demo{Proof} We have 
$[\eta]=\ddt(\exp\,t\eta) P_{\xi}\vert_{t=0}$, so
$i_2\circ i_1[\eta]=
\ddt \ga_{\xi}^{-1}(\exp\,t\eta) \ga_{\xi} 
\La^+\Gc \vert_{t=0} = [\la^{-i}\eta]$.
Similarly
$[\la^j\eta]=\ddt(\exp\,t\la^j \eta) 
\La_{\alg}^+\Gc\cap \ga_{\xi}(\La^+\Gc)\ga_{\xi}^{-1}
\vert_{t=0}$, so $i_2[\la^j\eta]$ is equal to
$\ddt \ga_{\xi}^{-1}(\exp\,t\la^j \eta) \ga_{\xi} 
\La^+\Gc \vert_{t=0} = [\la^{j-i}\eta]$.
\qed\enddemo 

The derivative $Du_{\xi}:T^{1,0}U_{\xi}\to T^{1,0}\Om_{\xi}$ 
corresponds to a map 
$$
\La_{\alg}^+\gc/
\La_{\alg}^+\gc\cap \Ad(\ga_{\xi})\La^+\gc
\to \gc/\p_{\xi}.
$$

\proclaim{Lemma 2.6}
The map $Du_{\xi}$ is given by
$$
[\la^{j}\eta]\mapsto
\cases
0 &\ \text{if}\ j>0\\
[\eta] &\ \text{if}\ j=0
\endcases
$$
where $\eta\in\g^{\xi}_i$.
\endproclaim

\demo{Proof} This follows from a calculation similar to that 
of Lemma 2.5, noting that $u_{\xi}$ is defined by 
evaluation at $\la=0$.
\qed\enddemo

We shall also need later the following mild generalization of
Proposition 2.4, describing the part of $U_{\xi}$ which lies 
over the \ll big cell\rr of $\Om_{\xi}$. The big cell of
$\Om_{\xi}$ means the subspace
$$
N'\cdot \ga_{\xi} \ \sub\ \Om_{\xi}
$$
where $N'$ is the \ll opposite\rr nilpotent subgroup to $N$,
i.e. the connected subgroup of $G$ with Lie algebra
$$
\n'={\tsize\bigoplus}_{\al\in\De^+}\gal.
$$
By a calculation similar to that of Proposition 2.4,
$N'\cdot \ga_{\xi}$ is diffeomorphic to the vector space
$\bigoplus_{\al\in\De^+,\,\al(\xi)\ne 0}\gal$, i.e. 
$\bigoplus_{i>0}\g^{\xi}_i$. 
The part of $U_{\xi}$ which lies 
over this \ll big cell\rr is
$$
u_{\xi}^{-1}(N'\cdot \ga_{\xi})=
\La_{N',\alg}^+\Gc\cdot\ga_{\xi}
$$
where $\La_{N',\alg}^+\Gc=\{\ga\in\La_{\alg}^+\Gc
\st \ga(0)\in N' \}$. Using the method of Proposition 2.4
again, we have:

\proclaim{Proposition 2.7}
$\La_{N',\alg}^+\Gc\cdot\ga_{\xi}=
\exp\,\u^0_{\xi}\cdot \ga_{\xi}\cong \exp\,\u^0_{\xi}
\cong \u^0_{\xi}$,
where $\u^0_{\xi}$ is the (finite dimensional) subalgebra 
$\bigoplus_{0\le i<r(\xi)}\la^i(\f^{\xi}_i)^\perp$ of
$\La^+_{\alg}\gc$.
\qed
\endproclaim

\newpage
\head
\S 3 The twistor construction
\endhead

\subheading{$S^1$-invariant extended solutions}

A fundamental role in the classification theory is played
by the following special extended solutions:

\no{\it Definition:} An  {\it $S^1$-invariant extended solution} 
is an extended solution $\Phi:M\to \Om G$ such that
$\Im\,\Phi\sub\Om_{\xi}$, for some $\xi\in I'$.

Let $\Phi:M\to\Om_{\xi}$ be an ($S^1$-invariant) extended 
solution. From the definition of $\Om_{\xi}$, the harmonic
map $\phi=\Phi(\ ,-1):M\to G$ factors through
$$
N_{\xi}=\{g\ga_{\xi}(-1)g^{-1}\st g\in G\},
$$
i.e. the conjugacy class of $\ga_{\xi}(-1)$. This is
a symmetric space; it is diffeomorphic to $G/C(\ga_{\xi}(-1))$, where
$C(g)$ denotes the centralizer of $g$. The inclusion of
$N_{\xi}$ in $G$ is known to be totally geodesic, with respect 
to the natural Riemannian metrics on $N_{\xi}$ and $G$ constructed from
a bi-invariant inner product on $\g$. It follows from this that a
harmonic map $M\to N_{\xi}$ is the same thing as a 
harmonic map $M\to G$ which factors through $N_{\xi}$. Hence,
$\phi$ is a harmonic map into $N_{\xi}$. Observe that
$\phi=\pi_{\xi}\circ\Phi$,  where the map 
$$
\pi_{\xi}:\Om_{\xi}\to N_{\xi}
$$ 
is given by $g\ga_{\xi}g^{-1}\mapsto g\ga_{\xi}(-1)g^{-1}$.

The extended
solution equation admits the following interpretation 
in this case. We define a subbundle
$$
H^{1,0}\Om_{\xi} = \Gc\times_{P_{\xi}}(\f^{\xi}_1/\f^{\xi}_0)
$$
of the holomorphic tangent bundle 
$T^{1,0}\Om_{\xi} \cong \Gc\times_{P_{\xi}}(\gc/\f^{\xi}_0)$,
and we say that a holomorphic map
$\Phi:M\to\Om_{\xi}$ is {\it super-horizontal} (with respect to
$\pi_{\xi}$) if and only
if $\Phi_z$ takes values in $H^{1,0}\Om_{\xi}$. Then we have:

\proclaim{Proposition 3.1}
A holomorphic map
$\Phi:M\to\Om_{\xi}$ is an extended solution if and
only if it is super-horizontal.
\endproclaim

\demo{Proof}
This follows immediately from Proposition 1.4 and Lemma 2.5.
\qed\enddemo

\no Hence, if $\Phi:M\to \Om_{\xi}$ is holomorphic and
super-horizontal, then $\pi_{\xi}\circ\Phi:M\to N_{\xi}$ 
is harmonic. In the terminology of \cite{Bu-Ra}, 
$\pi_{\xi}:\Om_{\xi}\to N_{\xi}$ is a {\it twistor fibration}.

If $\Phi=[\Psi\ga_{\xi}]$, where $\Psi:M\to\Gc$, then the
condition for $\Phi$ to be super-horizontal and holomorphic
may be written more explicitly as
$$
\align
\Im \Psi^{-1}\Psi_z &\sub \f^{\xi}_1\\
\Im \Psi^{-1}\Psi_{\bar z} &\sub \f^{\xi}_0.
\endalign
$$
This is a special case of Proposition 1.4.

The uniton number of an $S^1$-invariant extended solution
is given by:

\proclaim{Proposition 3.2}
Let $\Phi:M\to\Om_{\xi}$ be an extended solution.  Then
the uniton number of $\Phi$ is $r(\xi)$.
\endproclaim

\demo{Proof} It suffices to prove that $\Ad (g\ga_{\xi})$
is of the form $\sum_{\vert i\vert\le r(\xi)} \la^i T_i$,
with $T_{r(\xi)}\ne 0$, for any $g\in G$.  In fact, it suffices
to prove this when $g=e$. We know that
$\Ad\ga_{\xi}$ is given by multiplication by 
$\la^i$ on $\g^{\xi}_i$. Since 
$r(\xi)=\max\{i\st \g^{\xi}_i\ne 0\}$, the result follows.
\qed\enddemo

\subheading{Canonical twistor fibrations}

Harmonic maps into symmetric spaces which correspond to
$S^1$-invariant extended solutions have been studied
intensively.  It turns out that they arise from a distinguished
subset of the twistor fibrations $\pi_{\xi}$, the so called
canonical twistor fibrations.  These fibrations may be
described without reference to the loop group of $G$.
We shall review this theory, following \cite{Bu-Ra}. 

Let $\gc$ be a complex semisimple Lie algebra, with 
(compact) real form
$\g$, Cartan subalgebra $\tc$, and simple roots $\al_1,\dots,\al_l$.  
Any subset $\Cal P= \{\al_{i_1},\dots,\al_{i_k}\}$ 
of $\{\al_{1},\dots,\al_{l}\}$ defines a subalgebra $\p$ of
$\gc$, via the formula
$$
\p=\tc\oplus
(
{\tsize\bigoplus}_{\al\in\De^-}\gal
)
\oplus
(
{\tsize\bigoplus}_{\al\in\De'}\gal
)
$$
where   
$\De'=\{\al\in\De^+ \st \al=\sum_{i\notin \Cal P} n_i\al_i,
n_i\ge 0\}$.
It is well known that this gives a one to one correspondence
between subsets of the simple roots and (conjugacy classes
of) parabolic subalgebras of $\gc$.

Let $\xi_1,\dots,\xi_l\in\t$ be dual to $\al_1,\dots,\al_l$,
in the sense that $\al_i(\xi_j)=\psq \de_{ij}$.  Let $G$ be
the compact connected Lie group with {\it trivial centre}
corresponding to the Lie algebra $\g$.  The lattice
$\Z\xi_1\oplus\dots\oplus\Z\xi_l$ is the integer lattice $I$
of $G$.  To the subset  $\Cal P= \{\al_{i_1},\dots,\al_{i_k}\}$
we associate the element $\xi=\xi_{i_1}+\dots+\xi_{i_k}$ of
$I$. In the notation of \S 2, we then have 
$$
\p=\bigoplus_{i\le 0}\,\g^{\xi}_i=\p_{\xi}.
$$
Let $\h=\p\cap\g$, and let $H$ be the connected subgroup
of $G$ with Lie algebra $\h$. Evidently we have $\hc=\g^{\xi}_0$.
The natural inclusion $G/H\to \Gc/P$ is a diffeomorphism
(for dimensional reasons). 

We define a subalgebra $\kc$ of $\gc$ by
$$
\kc=\bigoplus_{i\,\text{even}}\,\g^{\xi}_i.
$$  Let
$\k=\kc\cap\g$, and let $K$ be the connected subgroup of
$G$ with Lie algebra $\k$. The homogeneous space
$G/K$ is a symmetric space.  Since $K\sub H$,
there is a  natural map $\pi:G/H\to G/K$. The map $\pi$
is called the {\it canonical twistor fibration} 
associated to $\Cal P$.

\proclaim{Proposition 3.3 } Let $\xi=\xi_{i_1}+\dots+\xi_{i_k}$,
as above. Then the map $\pi:G/H\to G/K$ coincides with the
map $\pi_{\xi}:\Om_{\xi}\to N_{\xi}$.
\qed
\endproclaim

\demo{Proof} We have already noted that $\p=\p_{\xi}$.
Since $N_{\xi}$ is the conjugacy class of $\ga_{\xi}(-1)$,
it suffices to prove that the Lie algebra of the
centralizer of $\ga_{\xi}(-1)$
is $\k$. This is elementary.
\qed
\enddemo

In view of this, we make the following definition:

\no{\it Definition:} Let $\xi\in I'$. We say that 
$\xi$ is {\it canonical} if $\xi=\xi_{i_1}+\dots+\xi_{i_k}$
for some $i_1,\dots,i_k$ (equivalently, if each simple root
takes the value $0$ or $\psq$ on $\xi$). We call
$\ga_{\xi}$ a {\it canonical geodesic}.

\no In general, if $\xi\in I'$, then we have 
$\xi=\sum_{j=1}^k n_{i_j}\xi_{i_j}$ for some 
$n_{i_1},\dots,n_{i_k}\in\N$. So $\xi$ is canonical
if and only if each $n_{i_j}$ is $1$.  We shall 
make use of the following technical result later on.

\proclaim{Lemma 3.4 \cite{\bf Bu-Ra}}
Let $r=r(\xi)=\max\{i\st\g^{\xi}_i\ne0\}$.

\no (1) If $\xi$ is canonical, then 
$\g^{\xi}_i\ne 0$ if and only if $-r\le i\le r$.

\no(2) Let $\xi=\sum_{j=1}^k n_{i_j}\xi_{i_j}$. We define
$\xi_{\can}$ by $\xi_{\can}=\sum_{j=1}^k \xi_{i_j}$. Then
$\g^{\xi}_0=\g^{\xi_{\can}}_0$ and  
$\f^{\xi}_0=\f^{\xi_{\can}}_0$.
\qed
\endproclaim

\newpage
\head
\S 4 Harmonic maps into Lie groups
\endhead

\subheading{Classification of extended solutions in terms
of canonical elements}

The basis for our analysis of extended solutions of
finite uniton number is the following simple observation:

\proclaim{Proposition 4.1}
Let $\Phi:M\to \Om_{\alg}^k G$ be an extended solution.
Then there exists some $\xi\in I'$, and some discrete subset
$D$ of $M$, such that $\Phi(M-D)\sub U_{\xi}$.
\endproclaim

\demo{Proof} This is a consequence of three facts:
(i) $\Phi$ is holomorphic,
(ii) $\dim_{\C}M=1$, and
(iii) the closures of the pieces of the
Bruhat decomposition give algebraic subvarieties of the 
(finite dimensional) complex projective
algebraic variety $\Om_{\alg}^k G$.
\qed\enddemo

Recall (from \S 2) that we have a holomorphic vector bundle 
$u_{\xi}:U_{\xi}\to\Om_{\xi}$.

\proclaim{Proposition 4.2}
If $\Phi:M-D\to U_{\xi}$ is an extended solution,
then $u_{\xi}\circ\Phi:M-D\to \Om_{\xi}$ is an extended solution.
\endproclaim

\demo{Proof}
Certainly $u_{\xi}\circ\Phi$ is holomorphic. By the formula
for $Du_{\xi}$ of Lemma 2.6,
$u_{\xi}\circ\Phi$ is an extended solution.
\qed\enddemo

\proclaim{Proposition 4.3}
If $\Phi:M-D\to U_{\xi}$ is an extended solution,
then the uniton number of $\Phi$ is equal to
the uniton number of $u_{\xi}\circ\Phi$, namely $r(\xi)$.
\endproclaim

\demo{Proof} This is similar to the proof of Proposition 3.2.
If $\ga\in\La_{\alg}^+\Gc$, we have  
$\ga\cdot\ga_{\xi}=\ga\ga_{\xi}\de$ for
some $\de\in\La_{\alg}^+\Gc$. In particular, 
we have $\ga(0),\de(0)\in\Gc$.
By the argument of Proposition 3.2, we see that
$\Ad(\ga\cdot\ga_{\xi})$  is of the form
$\sum_{i\ge -r(\xi)} \la^i T_i$, with $T_{-r(\xi)}\ne 0$.
Since $\Ad(\ga\cdot\ga_{\xi})\in O(\gc)$, we must
have $T_i=\bar T_{-i}$ for all $i$, so the result follows.
\qed\enddemo

Thus, to any extended solution $\Phi$ (of finite
uniton number) we may associate an $S^1$-invariant
extended solution $u_{\xi}\circ\Phi$, with the
same uniton number. Geometrically, $u_{\xi}\circ\Phi$
is obtained from $\Phi$ by applying the gradient flow of the
Morse-Bott function $E:\Om G\to \R$.  From the description
of this flow in \S 2, we have $u_{\xi}\circ\Phi=
\lim_{t\to\infty}\Phi^t$, where $\Phi^t$ is the extended
solution $e^{-t}\cdot\Phi$. This gives
a \ll classification\rr of extended solutions in terms of
$S^1$-invariant extended solutions, which is finer than
the classification by uniton number.

As a first step in using this classification, 
we note a useful reformulation of the
extended solution equation. For any (smooth) map
$\Phi:M-D\to U_{\xi}$, and any point $z_0\in M-D$,
we may write
$$
\Phi\vert_{M_0}=A\cdot\ga_{\xi},\quad
A:M_0\to \La_{\alg}^+\Gc
$$
where $M_0$ is some neighbourhood of $z_0$ in $M-D$.
In other words, we
represent $\Phi:M_0\to \Om G\cong \La\Gc/\La^+\Gc$
as $[A\ga_{\xi}]$, where $A\ga_{\xi}:M_0\to \La_{\alg}\Gc$.
By Proposition 1.4, the extended solution equation is 
equivalent to the conditions
$$
\align
\Im\,\la \Ad(\ga_{\xi}^{-1})A^{-1}A_z&\sub \La^+\gc \\
\Im\,\Ad(\ga_{\xi}^{-1})A^{-1}A_{\bar z}&\sub \La^+\gc.
\endalign
$$
Let us write
$$
A^{-1}A_z=\sum_{i\ge 0}\la^i A'_i,\quad
A^{-1}A_{\bar z}=\sum_{i\ge 0}\la^i 
A^{\prime\prime}_i.
$$
Since $\Ad(\ga_{\xi}^{-1})\eta=\la^{-i}\eta$
if $\eta\in\g^{\xi}_i$, we have:

\proclaim{Proposition 4.4}
The extended solution equation for 
$\Phi=[A\ga_{\xi}]:M_0\to U_{\xi}$ is
equivalent to the conditions
$$
\align
\Im\,A'_i&\sub \f^{\xi}_{i+1}\quad
\text{for}\ 0\le i\le r(\xi)-2\\
\Im\,A^{\prime\prime}_i&\sub \f^{\xi}_{i}\quad\ \ \,
\text{for}\ 0\le i\le r(\xi)-1.
\endalign
$$
(For higher values of $i$, these conditions are vacuous.)
\qed
\endproclaim

\subheading{Estimates of the minimal uniton number}

The possible values of the minimal uniton number of a
harmonic map are severely restricted by the next theorem.

\proclaim{Theorem 4.5}
Assume that $G$ is semisimple, with trivial centre.
Let $\Phi:M\to \Om_{\alg}^k G$ be an extended solution 
(of finite uniton number).
Then there exists some {\rm canonical} $\xi\in I'$,
some $\ga\in \Om_{\alg} G$,  and some discrete subset
$D$ of $M$, such that $\ga\Phi(M-D)\sub U_{\xi}$.
\endproclaim

\demo{Proof} As explained above, we may represent 
$\Phi\vert_{M-D}$ as $A\cdot\ga_{\xi}$ for some 
$\xi\in I'$, some 
discrete subset $D$ of $M$, and some $A:M-D\to\La_{\alg}^+\Gc$ 
satisfying the conditions of Proposition 4.4. Since
$G$ is semisimple, with trivial centre, we may write
$\xi=\sum_{j=1}^k n_{i_j}\xi_{i_j}$, with $n_{i_j}\ge 1$, where
$\xi_1,\dots,\xi_l$ are dual to the simple roots 
$\al_1,\dots,\al_l$ (as in \S 3). If $\xi$
happens to be canonical, i.e. $n_{i_j}=1$ for all $j$,
we are done. If not, then 
$\hat\xi=\xi-\xi_{\can}=\sum_{j=1}^k (n_{i_j}-1)\xi_{i_j}$ 
is a non-zero element of $I'$.

We claim that  
$\f^{\xi}_{i+1}\sub \f^{\hat\xi}_i\ \text{for}\ i\ge 0$,
and $\f^{\xi}_{0}\sub \f^{\hat\xi}_0$.
Assuming this for the moment, we deduce from Proposition 4.4
that $A$ satisfies the conditions
$$
\align
\Im\,A'_i&\sub \f^{\hat\xi}_{i}\\
\Im\,A^{\prime\prime}_i&\sub \f^{\hat\xi}_{i}
\endalign
$$
for $0\le i\le r(\xi)-1$.
These conditions say that $A\cdot \ga_{\hat\xi}$ is 
both holomorphic and anti-holomorphic, hence constant (in $z$).
In other words, it defines an element of 
$\Om G$, which we shall call $\ga^{-1}$. 
Thus, $\ga^{-1}=A\cdot \ga_{\hat\xi}=A\ga_{\hat\xi}B$,
for some map $B:M-D\to\La_{\alg}^+\Gc$. We then have
$\ga\Phi\cdot\ga_{\xi}=
B^{-1}\ga_{\hat\xi}^{-1}A^{-1}A\cdot\ga_{\xi}
=B^{-1}\cdot \ga_{\xi_{\can}}$,
which is the desired conclusion.

It remains to prove the above claim. Since $\Ad\ga_{\xi}$, 
$\Ad\ga_{\xi_{\can}}$ and 
$\Ad\ga_{\hat\xi}$ are simultaneously diagonalized on the
root space $\gal$, it follows that
$$
\g^{\xi}_k={\tsize\bigoplus}_{0\le i\le k}
\left(
\g^{\hat\xi}_i \cap \g^{\xi_{\can}}_{k-i}
\right).
$$
{}From this and Lemma 3.4 (2),
$\f^{\xi}_i\sub \f^{\hat\xi}_i$ for all $i\ge 0$.
To show that $\f^{\xi}_{i+1}\sub \f^{\hat\xi}_i$
for $i\ge 0$, it suffices to show that
$\g^{\xi}_{i+1}\cap(\g^{\hat\xi}_{i+1}\cap\g^{\xi_{\can}}_0)=0$.
{}From Lemma 3.4 (2) again, we have 
$\g^{\xi}_0=\g^{\xi_{\can}}_0$, so
$\g^{\xi}_{i+1}\cap\g^{\xi_{\can}}_0=
\g^{\xi}_{i+1}\cap\g^{\xi}_0=0$.
\qed\enddemo

The key point of this theorem is that we 
reduce to the situation of a {\it canonical} $\xi$. 
(This reduction is analogous to the concept of 
\ll normalization\rr in \cite{Uh} and \cite{Se},
and to the geometrical idea of reducing to \ll full\rr harmonic
maps in earlier works on this subject.) 
As a consequence, we obtain estimates
for the minimal uniton number of a harmonic map:

\proclaim{Corollary 4.6}

\no(1) Let $G$ be any compact simple Lie group.
Let $\phi:M\to G$ be a harmonic map. Assume that $\phi$
arises from an extended solution $\Phi:M\to\Om G$ of finite
uniton number. (If $M=S^2$, this assumption holds automatically.) 
Then the minimal uniton number of $\phi$ is
not greater than $r(G)=\sum_{i=1}^l n_i$, where 
$\al=\sum_{i=1}^l n_i\al_i$ is the expression for the highest
root of $G$ in terms of the simple roots $\al_1,\dots,\al_l$.

\no(2) Let $G$ be any compact Lie group, and let the
simple factors in the universal cover of $G$ be $G_1,\dots,G_t$.
Let $\phi$ be as in (1). Then the minimal uniton number of 
$\phi$ is not greater than $\max\{r(G_1),\dots,r(G_t)\}$.
\endproclaim

\demo{Proof} (1) Since the minimal uniton number depends only on
$\Ad\phi:M\to\Ad G$, it suffices to prove the statement when
$G$ has trivial centre.
Using the notation of \S 2 and \S 3, we have
$$
\align
\max&\{r(\xi)\st \xi\ \text{canonical}\}\\
&=
\max\left\{ \al(\sum_{i\in\Cal P}\xi_i)/\psq 
\ \right|\  \left.
\Cal P\sub \{\al_1,\dots,\al_l\},\ 
\al=\sum_{i=1}^l m_i\al_i\in\De^+\right\}\\ 
&=(\sum_{i=1}^l n_i\al_i)(\sum_{i=1}^l\xi_i)/\psq
=\sum_{i=1}^l n_i.
\endalign
$$
By Proposition 3.2 (and the above theorem), this is an upper
bound for the minimal uniton number.
(2) As in part (1), it suffices to prove the statement for
$\Ad G_1\times\dots\times\Ad G_t$. But the statement in this
case follows immediately from (1).
\qed\enddemo

We list below the values of $r(G)$ for (representatives of local
isomorphism classes of) the compact simple Lie groups.

\bigpagebreak

\settabs 5\columns

\+  & $G$ & $r(G)$ & $d(G)-1$ &\cr

\+ & -------- & -------- & -------- &\cr

\+ & $SU_n$ & $n-1$ & $n-1$ &\cr

\+ & $SO_{2n+1}$ & $2n-1$ & $2n$  &\cr

\+ & $Sp_n$ & $2n-1$ & $2n-1$ &\cr

\+ & $SO_{2n}$ & $2n-3$ & $2n-1$ &\cr

\+ & $G_2$ & $5$ & $6$ &\cr

\+ & $F_4$ & $11$ & $25$ &\cr

\+ & $E_6$ & $11$ & $26$ &\cr

\+ & $E_7$ & $17$ & $55$ &\cr

\+ & $E_8$ & $29$ & $247$ &\cr

The fact that $r(SU_n)=n-1$ was proved in \cite{Uh} and \cite{Se}.
If $\theta:G\to SU_n$ is a faithful representation, then
$n-1$ is an upper bound for the minimal uniton number of a
harmonic map into $G$. By taking $n=d(G)$, the dimension of
the smallest faithful representation of $G$,
we obtain the upper bounds $d(G)-1$ in the third column. 
Our method shows that these upper bounds can be
sharpened, especially for the exceptional groups.
The question arises as to whether these upper bounds are
optimal; it turns out that they are:

\proclaim{Proposition 4.7} For any compact simple Lie
group $G$, there exists a harmonic map $\phi:S^2\to G$ 
such that $r(\phi)=r(G)$.
\endproclaim

\demo{Proof} It suffices to construct $\phi$ when 
$G$ has trivial centre. Let $\xi_1,\dots,\xi_l$ be as above,
and let $\xi=\xi_1+\dots+\xi_l$. We shall construct a
holomorphic super-horizontal map $\Phi:S^2\to \Om_{\xi}
\cong \Gc/P_{\xi}$ of minimal uniton number $r(G)$.

For each $i=1,\dots,l$, let $X_i$ be a non-zero vector in 
$\g_{\al_i}$, and let $X=X_1+\dots+X_l \in \g^{\xi}_1$.
By \cite{Ko}, there exists some $Y \in \g^{\xi}_{-1}$
such that $\Span\{\xi,X,Y\}$ is a subalgebra isomorphic
to $\frak{sl}_2\C$.  Let $S_{\xi}$ be the corresponding
subgroup (isomorphic to $SL_2\C$). The induced map
$\Phi:S_{\xi}/S_{\xi}\cap P_{\xi} \to \Gc/P_{\xi}$ is
holomorphic and super-horizontal (see \S 3). To prove
the proposition, we will show that $r(\ga\Phi)\ge r(\xi)$
for all $\ga\in\Om G$.

Let us assume that $r(\ga\Phi) < r(\xi)$, for some $\ga$.
It is a property of $SL_2\C$ that there exists some
$g\in S_{\xi}$ such that $\Ad(g)\xi=-\xi$, hence for this
$g$ we have $\Phi([g])=\Phi([e])^{-1}$.  It follows
that $\Phi([e])^{2}=\Phi([g])^{-1}\Phi([e])=
(\ga\Phi)([g])^{-1}(\ga\Phi)([e])$. If $r(\ga\Phi) < r(\xi)$,
this is a contradiction.
\qed
\enddemo

As a consequence of Theorem 4.5, we can elucidate the structure
of harmonic maps $\phi$ with $r(\phi)=0$, $1$, or $2$. It is
clear that $r(\phi)=0$ if and only if $\phi$ is constant.
The situation for $r(\phi)=1$ is also easy (and well known):

\proclaim{Corollary 4.8}
Let $\phi:M\to G$ be harmonic. Then $r(\phi)=1$ if and only if
$\phi$ is (up to left translation by an element of $G$)
the composition of a holomorphic map
$M\to G/K$, for some Hermitian symmetric space $G/K$, with
a totally geodesic embedding of $G/K$ in $G$.
\endproclaim

\demo{Proof} Suppose $r(\phi)=1$. For $r(\xi)=1$, the formula
for the rank of the bundle $U_{\xi}\to\Om_{\xi}$ shows that
$U_{\xi}=\Om_{\xi}$ (i.e. the Morse index is zero). Moreover,
from \S 3 we see that the twistor fibration $\Om_{\xi}\to N_{\xi}$
is the identity map. The symmetric space $N_{\xi}$ is Hermitian
($\ad\xi$ provides an invariant complex structure), so $\phi$
is of the stated form.  The converse statement follows 
directly from the existence of a canonical (trivial) twistor 
fibration for any Hermitian symmetric space (see \cite{Bu-Ra}).
\qed\enddemo

For $r(\phi)=2$, we have:

\proclaim{Corollary 4.9}
Let $\xi\in I'$ with $r(\xi)\le 2$.
A holomorphic map $\Phi:M-D\to U_{\xi}$ is an extended solution
if and only if $u_{\xi}\circ\Phi:M-D\to \Om_{\xi}$ is
super-horizontal. 
\endproclaim

\demo{Proof} The condition for $\Phi$ to be an extended
solution is $\Im A'_0\sub \f^{\xi}_1$. This is the condition 
for $\Phi$ to be super-horizontal. (The conditions
on $A^{\prime\prime}$ are satisfied because $\Phi$ is
holomorphic.)
\qed\enddemo

\no This shows that harmonic maps with $r(\phi)=2$ correspond to
pairs $(\Phi,\sigma)$, where $\Phi$ is a holomorphic 
super-horizontal map into $\Om_{\xi}$ and $\sigma$ is a
meromorphic section of the holomorphic vector bundle
$\Phi^\ast U_{\xi}$. For example, all harmonic maps
$S^2\to SU_3$ or $U_3$ are of this form.

The situation for $r(\phi)\ge 3$ is inevitably more complicated.
However, we shall see that a rather explicit description
is possible even in this case.

\subheading{Weierstrass formulae for extended solutions}

Further information may be obtained by making use of 
the (Zariski) open subset $u_{\xi}^{-1}(N'\cdot\ga_{\xi})
=\exp\,\u^0_{\xi}\cdot\ga_{\xi}$ of $U_{\xi}$, from
Proposition 2.7. Since an extended solution
$\Phi:M-D\to U_{\xi}$ is holomorphic, and the complement
of $u_{\xi}^{-1}(N'\cdot\ga_{\xi})$ in $U_{\xi}$
is a proper algebraic subvariety, there exists a discrete
subset $D'$ of $M$, with $D\sub D'\sub M$, such that
$\Phi(M-D')\sub \exp\,\u^0_{\xi}\cdot\ga_{\xi}$. 
(After left translation by a constant loop if necessary,
the image of $\Phi$ will not be entirely
contained in the complement of 
$u_{\xi}^{-1}(N'\cdot\ga_{\xi})$.) We can therefore write
$$
\Phi\vert_{M-D'}=\exp\,C\cdot\ga_{\xi}
$$
where $C:M-D'\to\u^0_{\xi}$ is a (vector-valued)
{\it holomorphic function}. The map $\Psi=(\exp\,C)\ga_{\xi}$
is a complex extended solution in the sense of \S 1, and what 
we have just said constitutes a proof of Theorem 1.5,
for extended solutions of finite uniton number. Since 
$\Om^k_{\alg}G$ is
a projective algebraic variety, $C$ is a meromorphic 
function from $M$ to $\u^0_{\xi}$.

Conversely, if $C:M\to\u^0_{\xi}$ is meromorphic, the
condition (from Proposition 4.4) for 
$\Phi=\exp\,C\cdot\ga_{\xi}$ to be an extended solution is
$$
\Im A'_i\sub \f^{\xi}_{i+1},\quad 0\le i\le r(\xi)-2
$$
where $(\exp\,C)^{-1}(\exp\,C)_z
=\sum_{0\le i\le r(\xi)-1}\la^i A'_i$.
We shall investigate this condition more closely.
To do so, let us write $r=r(\xi)$ and
$$
C= C_0 + \la C_1 + \dots + \la^{r-1}C_{r-1},\quad
C_i = c^{i+1}_i + c^{i+2}_i + \dots + c^r_i
$$
where each function $c^j_i:M\to \g^{\xi}_j$ is meromorphic. We
shall make essential use of the fact that
$[\la^i \g^{\xi}_{i+k},\la^j \g^{\xi}_{j+l}]\sub
\la^{i+j} \g^{\xi}_{i+j+k+l}$.

By the well known formula for the derivative of
the exponential map (\cite{He}, Chapter 2, Theorem 1.7), we have
$$
\align
(\exp\,C)^{-1}(\exp\,C)_z&=\tfrac{I-e^{-\ad C}}{\ad C}C_z\\
&=C_z-\tfrac{1}{2!}(\ad C)C_z+\tfrac{1}{3!}(\ad C)^2C_z-
\tfrac{1}{4!}(\ad C)^3C_z+\dots
\endalign
$$
The condition for an extended solution is that the
coefficient of $\la^i$ in this expression should
have zero component in each of 
$\g^{\xi}_{i+2},\dots,\g^{\xi}_{r}$.

For $i=0$, this means that
$$
(C_0)_z-\tfrac{1}{2!}(\ad C_0)(C_0)_z+\tfrac{1}{3!}
(\ad C_0)^2(C_0)_z-\tfrac{1}{4!}(\ad C_0)^3(C_0)_z+\dots
$$
should have zero component in each of 
$\g^{\xi}_{2},\dots,\g^{\xi}_{r}$.  We obtain 
equations (for $j=2,\dots,r$) of the form 
\ll$(c^{j}_0)_z=$ terms involving
$c^{1}_0,\dots,c^{j-1}_0$ and their $z$-derivatives\rrr.
There is no condition at all on $c^{1}_0$, which
may be taken to be any $\g^{\xi}_1$-valued meromorphic 
function on $M$. Each of $c^{2}_0,\dots,c^{r}_0$ may then
be determined (locally) by integration.

For the coefficient of $\la^i$, when $i>0$, 
we obtain a similar system of
equations for $C_i$, i.e. for $c^{i+2}_i,\dots,c^r_i$.
For each $j=i+2,\dots,r$, we have an equation of the form 
\ll$(c^{j}_i)_z=$ terms involving
$c^{k}_l$ for $k<i+2,l<i$ and their $z$-derivatives\rrr.
Therefore we may choose $c^{i+1}_i$ to be any 
$\g^{\xi}_{i+1}$-valued meromorphic function on $M$,
after which $c^{i+2}_i,\dots,c^{r}_i$ may be
determined (locally) by integration.

Thus, any choice of meromorphic functions 
$c^1_0,c^2_1,\dots,c^r_{r-1}$ gives rise --- locally ---
to an extended solution of uniton number (at most) $r$.
It is not guaranteed that the extended solution is defined on the
whole of $M$.  To determine when this happens is a more delicate
matter (see \cite{Br1} for an example). Nevertheless, we
can conclude at least that all harmonic maps of $S^2$ are
\ll algebraic\rr in the following sense:

\proclaim{Theorem 4.10}
Every harmonic map $S^2\to G$ arises from an extended
solution which may be obtained explicitly by choosing a finite number
of rational functions and then performing a finite
number of algebraic operations and integrations.
\qed
\endproclaim

\no Weierstrass was the first to give a procedure of this
type, in the context of minimal surfaces in $\R^3$. 
Bryant (\cite{Br1}) used the above method to give
explicit constructions of harmonic maps from Riemann 
surfaces into $S^4$, from meromorphic functions.
In later work on the twistor construction (\cite{Br2}),
he gave essentially the above method for holomorphic
super-horizontal maps into $\Om_{\xi}$. For harmonic maps from
$S^2$ into $U_n$, Wood (\cite{Wo2}) proved a version of
Theorem 4.10,  by making a detailed analysis of Uhlenbeck's results.

\subheading{Example: harmonic maps $S^2\to U_n$}

With standard conventions for the maximal torus, integer lattice,
and fundamental Weyl chamber, the canonical geodesics
$\ga_{\xi}:S^1\to U_n$ are the homomorphisms of the form
$\ga_{k_1,k_2,\dots,k_n}(\la)=\diag(\la^{k_1},\la^{k_2},
\dots \la^{k_n})$, with $k_1\ge k_2\ge \dots \ge k_n$, and
$\vert k_i-k_{i-1}\vert=0$ or $1$ for all $i$. We may assume
that $k_n=0$; this corresponds to fixing representatives
of the canonical geodesics for $PU_n=U_n/Z(U_n)$, under the
surjection $\pi_1 U_n\cong\Z\to\Z/(n-1)\Z \cong \pi_1 PU_n$.
There are $2^{n-1}$ canonical geodesics of this type.

Let us consider the Weierstrass representation for the
most complicated type of harmonic map, namely where the
canonical geodesic is
$\ga_{\xi}(\la)=(\la^{n-1},\la^{n-2},\dots,1)$. (The minimal
uniton number takes its greatest possible value here, $n-1$.)
We have $\gc=\g^{\xi}_{n-1}\oplus\dots\oplus\g^{\xi}_{0}\oplus
\dots\oplus\g^{\xi}_{n-1}$, where $\g^{\xi}_{s}$ is the set of
$n\times n$ complex matrices $(x_{ij})$ with $x_{ij}=0$
unless $j-i=s$.

The nature of the equations for 
$C=C_0+\la C_1 +\dots + \la^{n-1}C_{n-1}$
will be more transparent if we write them out for $n=4$. In
this case, $C=C_0 + \la C_1 + \la^2 C_2$, and this is of the form
$$
C=
\pmatrix
{0} & a_1 & b_1 & c_1\\
{0} & {0} & a_2 & b_2\\
{0} & {0} & {0} & a_3\\
{0} & {0} & {0} & {0}
\endpmatrix
+\la 
\pmatrix
{0} & {0} & d_1 & e_1\\
{0} & {0} & {0} & d_2\\
{0} & {0} & {0} & {0}\\
{0} & {0} & {0} & {0}
\endpmatrix
+\la^2
\pmatrix
{0} & {0} & {0} & f_1\\
{0} & {0} & {0} & {0}\\
{0} & {0} & {0} & {0}\\
{0} & {0} & {0} & {0}
\endpmatrix.
$$
There are no conditions on $c^1_0=(a_1,a_2,a_3)$, 
$c^2_1=(d_1,d_2)$ and $c^3_2=f_1$. The conditions on the 
remaining components of $C_0$ and $C_1$ are as follows:
$$\gather
(C_0)_z-\frac12[C_0,(C_0)_z]+\frac16[C_0,[C_0,(C_0)_z]]
\ \text{has zero component in}\ \g^{\xi}_2,\g^{\xi}_3\\
(C_1)_z-\frac12([C_0,(C_1)_z]+[C_1,(C_0)_z])
\ \text{has zero component in}\ \g^{\xi}_3.
\endgather
$$

The first equation is simply the condition that
$\exp\,C_0\cdot\diag(\la^3,\la^2,\la,1)$ be super-horizontal.
The conjugacy class $\Om$ of the
canonical geodesic $\diag(\la^{n-1},\la^{n-2},\dots,1)$
is the full flag manifold $F_{1,2,\dots,n-1}(\C^n)$,
and it is easy to see that a super-horizontal holomorphic map 
$\Phi:S^2\to F_{1,2,\dots,n-1}(\C^n)$
may be represented on the complement of a finite set $D$
by a \ll holomorphic frame\rr of the form
$$
A=
\pmatrix
\vert &    & \vert & \vert\\
f^{(n-1)} & \dots & f' &  f  \\
\vert &    & \vert & \vert
\endpmatrix
$$
where $f:S^2-D\to\C^n$ is holomorphic. By our assumptions,
$$
A\cdot \diag(\la^3,\la^2,\la,1)=
\exp\,C_0\cdot\diag(\la^3,\la^2,\la,1).
$$
{}From this it is easy to verify that $\exp\,C_0$ must be of the form
$$
\exp\, C_0=
\pmatrix
1 & \de & \al'/\ga' & \al \\
{0} & 1 & \be'/\ga' & \be\\
{0} & {0} & 1 & \ga\\
{0} & {0} & {0} & 1
\endpmatrix,\quad
\de=(\al'/\ga')'/(\be'/\ga')'
$$
where $\al,\be,\ga$ are rational functions (and prime
denotes $z$-derivative). Conversely, for any $\al,\be,\ga$ ,
the above formula gives a solution $C_0$ of the first equation.

The second equation reduces to the differential equation
$$
e'_1=\frac12(a_1 d'_2 - a'_1 d_2 + d_1 a'_3 -d'_1 a_3),
$$
from which $e_1$ may be determined by integration.

Thus, any harmonic map $S^2\to U_4$ of this type corresponds
to six rational functions.

For $n=3$ there are no differential equations to solve
(beyond $C_0$),
as predicted by Corollary 4.9.  The most general harmonic
map with $r(\phi)=2$ arises from an extended solution
$\Phi=\exp\,C\cdot \diag(\la^2,\la,1)$, where
$$
C=
\pmatrix
1 & \al'/\be' & \al\\
{0} & 1 & \be\\
{0}  &  {0} & 1
\endpmatrix
+\la
\pmatrix
{0} & {0} & \ga\\
{0} & {0} & {0}\\
{0}  &  {0} & {0}
\endpmatrix
$$
and $\al,\be,\ga$ are arbitrary rational functions.
This description is equivalent to a description of
harmonic maps $S^2\to U_3$ given in \cite{Wo2}, 
although the latter 
was expressed somewhat differently, in terms of \ll uniton
factorizations\rrr. The role of such factorizations
will be considered next.
 
\subheading{Transforms and factorizations}

We have seen that the Bruhat decomposition of $\Om G$ provides
a straightforward approach to describing harmonic maps of
finite uniton number, particularly when the domain is $S^2$.
We shall indicate briefly how the existing theory of such 
maps can be understood from our point of view.

Initial work on harmonic maps from Riemann surfaces to
symmetric spaces was based on the idea of \ll B\"acklund
transformations\rrr, whereby a given harmonic map $\phi_1$ 
is transformed
into a new harmonic map $\phi_2$. Such transformations appear
in our theory in the following way. Consider an extended solution
$$
\Phi=A\cdot \ga_{\xi}:M-D\to U_{\xi},
\quad \xi=\xi_{i_1}+\dots+\xi_{i_k}.
$$
For any non-empty subset $J=\{j_1,\dots,j_m\}$ of 
$\{i_1,\dots,i_k\}$, we write 
$\xi_J=\xi_{j_1}+\dots+\xi_{j_m}$.

\proclaim{Proposition 4.11}
If $A\cdot\ga_{\xi}$ is an extended solution, then so
is $A\cdot\ga_{\xi_J}$.
\endproclaim

\demo{Proof}
It is obvious that
$\f^{\xi_J}_i\sub \f^{\xi}_i$ for all $i\ge 0$, so the assertion
follows from Proposition 4.4.
\qed
\enddemo

\no This gives $2^{k-1}$ transforms of our original extended 
solution.  For example, let us take $\Phi=A\cdot\ga_{\xi}$, where
$$
A=
\pmatrix
\vert &    & \vert & \vert\\
f^{(n-1)} & \dots & f' &  f  \\
\vert &    & \vert & \vert
\endpmatrix,\quad
\ga_{\xi}=
\diag(\la^{n-1},\la^{n-2},\dots,1)
$$
in the case $G=U_n$.  Then we obtain a new extended solution
$A\cdot \ga_{\xi_J}$ by taking
$\ga_{\xi_J}=\diag(\la^2,\dots,\la^2,\la,1,\dots,1)$,
where $\la$ appears as the $(n-i)$-th
diagonal element in $\ga_{\xi_J}$. In fact, $A\cdot \ga_{\xi_J}$
represents the \ll $i$-th harmonic transform\rr of the
holomorphic map $[f]:S^2\to\C P^{n-1}$, namely the map
$\Span\{f,f^\prime,\dots,f^{(i-1)}\}^\perp\cap
\Span\{f,f^\prime,\dots,f^{(i)}\}$.

In \cite{Uh}, Uhlenbeck considered harmonic maps into
complex Grassmannians as examples of harmonic maps into $U_n$, 
and introduced a transform procedure for such maps. 
This was defined as follows: if $\Phi_1,\Phi_2$ 
are extended solutions
corresponding to $\phi_1,\phi_2$, then $\Phi_2$ is obtained
from $\Phi_1$ by multiplication: $\Phi_2=U\Phi_1$, where
$U$ is a certain map into a Grassmannian (necessarily 
a solution to a first order system analogous 
to the Cauchy-Riemann equations).
The new extended solution $\Phi_2$ is said to be obtained 
from the old extended solution $\Phi_1$ by
\ll adding a uniton\rrr. To show that any harmonic map
can be constructed by transforming a constant map finitely
many times is equivalent to showing that any extended
solution can be factored in the form $\Phi=U_1\dots U_r$.
In \cite{Uh} and \cite{Se}, such a factorization was given
for harmonic maps from $S^2$ to $G=U_n$ or $G/K=Gr_k(\C^n)$.
This was subsequently generalized in \cite{Bu-Ra} to the
case of a group $G$ of \ll type $H$\rrr, i.e. a group
whose universal cover contains no factors locally isomorphic
to $G_2$, $F_4$, or $E_8$.

The factorization amounts to a reorganization of the earlier
B\"acklund transformation approach, and can be obtained from
our approach in a similar way. (We shall give the factorization
theorem for $U_n$ here, postponing the Grassmannian case to
the next section.) Namely, given $A\cdot \ga_{\xi}$,
we consider the sequence 
$$
\xi_{i_1},\ \xi_{i_1}+\xi_{i_2},\ \dots,
\ \xi_{i_1}+\dots+\xi_{i_k}.
$$
{}From this we obtain a sequence of extended solutions
$$
\Phi_1=A\cdot\ga_{i_1},\ \Phi_2=A\cdot\ga_{i_1}\ga_{i_2},\ 
\dots,\ \Phi_k=A\cdot\ga_{i_1}\ga_{i_2}\dots\ga_{i_k},
$$
where, to simplify notation, we write $\ga_{i}$ for 
$\ga_{\xi_i}$. When $G=U_n$, we claim that the factorization
$$
\Phi=\Phi_1(\Phi_1^{-1}\Phi_2)(\Phi_2^{-1}\Phi_3)\dots
(\Phi_{k-1}^{-1}\Phi_{k})
$$
is a \ll uniton factorization\rr in the sense of
\cite{Uh} or \cite{Se}, i.e. that 
$\Phi_j^{-1}\Phi_{j+1}$ is a map into a
Grassmannian.  To see this, we write
$\Phi_{j}=A\cdot\ga_{i_1}\dots\ga_{i_{j}}=
A\ga_{i_1}\dots\ga_{i_{j}} B_{j}$
for some $\La^+\Gc$-valued maps $B_j$. Then
$\Phi_j^{-1}\Phi_{j+1}=B_j^{-1}\ga_{i_{j+1}}B_{i_{j+1}}
=B_j^{-1}\cdot\ga_{i_{j+1}}$,
which is a map into $\La^+\Gc\cdot\ga_{i_{j+1}}$.
For $G=U_n$ (see Lemma 4.12 below), it turns out that
$\La^+\Gc\cdot\ga_{i_{j+1}}=\Gc\cdot\ga_{i_{j+1}}
=Gr_k(\C^n)$
for some $k$, so this completes the proof of our claim.

What we have shown is that a factorization of the extended solution
$\Phi=A\cdot\ga_{\xi}$ results from a choice of factorization
of the canonical geodesic $\ga_{\xi}$.  For example, when
$G=U_n$ and $\ga_{\xi}=\diag(\la^{n-1},\la^{n-2},\dots,1)$,
we may choose the factorization
$$
\diag(\la^{n-1},\la^{n-2},\dots,1)=
\diag(\la,\dots,\la,1)\ 
\diag(\la,\dots,\la,1,1)\ \ \dots\ \ 
\diag(\la,1,\dots,1).
$$ 

In order to generalize this to the case of other Lie groups,
we need:

\proclaim{Lemma 4.12}
Assume that $G$ is simple. Assume that $\al_i$ is a
miniscule simple root, i.e. that the coefficient of
$\al_i$ in the highest root of $G$ is one. Then
$\ga_{i}$ has Morse index zero (as a critical point
of the energy functional), and the critical
manifold $\Gc\cdot\ga_i$ is a Hermitian symmetric space.
\endproclaim

\demo{Proof} The assumption on $\al_i$ means that $r(\xi_i)=1$,
in our earlier notation. The conclusion follows as in the
proof of Corollary 4.8.
\qed\enddemo

\no When $G=SU_n$, all simple roots are miniscule, and
all Hermitian symmetric spaces $G/K$ are Grassmannians;
this is the fact that we used earlier.  For general $G$,
not all simple roots are miniscule --- in fact, the groups
$G_2$, $F_4$ and $E_8$ have no miniscule simple roots
at all. In the general case, the argument above gives:

\proclaim{Proposition 4.13}
Any extended solution $\Phi:M\to \Om_{\alg}G$ admits a
factorization $\Phi=\Phi_{1}\Phi_{2}\dots\Phi_{k}$, where
each $\Phi_{i}$ is a map into the (closure of the) unstable
manifold $U_{\xi_j}$ corresponding to a simple root $\al_j$.
Each sub-product $\Phi_{1}\Phi_{2}\dots\Phi_{i}$ is an
extended solution. 
\qed
\endproclaim

\no A refinement of such a factorization
into \ll linear factors\rr was given in \cite{Bu-Ra}, for
maps of $S^2$ into groups of type H, by exploiting further the
global geometry of $S^2$.

\newpage
\head
\S 5 Harmonic maps into symmetric spaces
\endhead

Our theory for Lie groups may be modified to deal with
symmetric spaces, if we replace loop groups by 
\ll twisted\rr loop groups. As in the case of Lie groups, we 
obtain short proofs of the known results,
as well as new results. We shall only consider {\it inner} symmetric
spaces, however. An inner symmetric space is a
homogeneous space of the form $G/K$ where
$C(g)_0\sub K\sub C(g)$ for some $g\in G$, where  
$C(g)_0$ is the identity component
of $C(g)$ (and $C(g)=\{x\in G\st xg=gx\}$). We shall use
the following convenient characterization up to local
isomorphism of such spaces.

\proclaim{Proposition 5.1}

\no(1) Let $G$ be a compact (connected) Lie group. 
Then each component of
$\sqrt e=\{g\in G\st g^2=e\}$ is a compact inner symmetric
space.  

\no(2) Conversely, any compact (connected) inner symmetric space
may be immersed in such a Lie group $G$ as a 
connected component of $\sqrt e$.  Moreover, it may be assumed
that $G$ has trivial centre.
\endproclaim

\demo{Proof} If $x\in\sqrt e$, then the conjugacy class of $x$
is an inner symmetric space, isomorphic to $G/C(g)$.
There are at most a finite number of such conjugacy classes
(as one sees by considering the intersection of $\sqrt e$
and a maximal torus of $G$). Hence $\sqrt e$ consists of
a finite number of conjugacy classes, each of which is a
connected component.  (2) This is a consequence of
Proposition 4.5 of \cite{Bu-Ra}.  The assertion concerning
the centre follows from the fact that the centre of $G$
is contained in the identity component
of $C(x)$, for any $x\in G$.
\qed\enddemo

\no For example, if $G=U_n$, then the connected components
of $\sqrt e$ are the symmetric spaces $Gr_k(\C^n)$ for
$k=0,1,\dots,n$.

It is well known that the embedding of each component of
$\sqrt e$ in $G$, in Proposition 5.1, is totally geodesic.  
Harmonic maps into
inner symmetric spaces may therefore be viewed as special
harmonic maps into $G$. As in \cite{Uh} and \cite{Se}, we
may characterize the corresponding special extended solutions
in terms of the involution
$$
T:\Om G\to\Om G,\quad T(\ga)(\la)=\ga(-\la)\ga(-1)^{-1}.
$$
We write
$$
(\Om G)_T=\{\ga\in\Om G\st T(\ga)=\ga\}
$$
for the fixed set of $T$. If $\Phi:M\to\Om G$ is an extended
solution, it is easy to verify that $T(\Phi)$ is also an
extended solution.

\proclaim{Proposition 5.2} 

\no(1) Let $\Phi:M\to (\Om G)_T$ be
an extended solution. Then $\phi(z)=\Phi(z,-1)$ defines a 
harmonic map from $M$ into (a connected component of) $\sqrt e$.

\no(2) Let $\phi:M\to\sqrt e$ be harmonic. Assume that there
exists an extended solution $\tilde\Phi:M\to\Om G$ such
that $\phi(z)=\Phi(z,-1)$. Then there exists an
an extended solution $\Phi:M\to(\Om G)_T$ such
that $\phi(z)=\Phi(z,-1)$.
\endproclaim

\demo{Proof} Statement (1) is obvious. For statement (2), 
we use the fact (from \S 1) that the set of extended solutions 
corresponding to $\phi$ is
$\{\ga\tilde\Phi\st \ga\in\Om G, \ga(-1)=e\}$.
We must find a $\ga$ such that $\Phi=\ga\tilde\Phi$ takes
values in $(\Om G)_T$. To do this, choose some $z_0\in M$,
and choose some $\xi\in\g$ with $\exp 2\pi\xi=e$ and 
$\phi(z_0)=\exp\pi\xi$. (This is possible as $\exp$ is
surjective and $\phi(z_0)^2=e$.) Then let 
$\ga=\ga_{\xi}\tilde\Phi(z_0)^{-1}$, so that $\Phi(z_0)=\ga_{\xi}$.
We claim that $\Phi$ takes values in $(\Om G)_T$. Both
$\Phi$ and $T(\Phi)$ are extended solutions corresponding 
to $\phi$, so we must have $\Phi=\de T(\Phi)$ for some
$\de\in\Om G$. But $\Phi(z_0)=T(\Phi)(z_0)=\ga_{\xi}$, so $\de=e$,
and $\Phi=T(\Phi)$, as required.
\qed\enddemo

The involution $T$ of $\Om G$ is holomorphic, as it is
induced by the complex involution $T:\La\Gc\to\La\Gc$
defined by $T(\ga)(\la)=\ga(-\la)$. The action of
$\C^\ast_{{\sssize\ge} 1}$ on $\Om G$ commutes with the
action of $T$ (since $-1\in \C^\ast_{{\sssize\ge} 1}$).
Hence, the Bruhat decomposition of $\Om_{\alg}G$ (described
in \S 2) induces a natural decomposition of $(\Om_{\alg}G)_T$.
For each $\xi\in I'$ (the intersection of the integer lattice
with a fundamental Weyl chamber) we obtain a holomorphic
vector bundle $(U_{\xi})_T\to (\Om_{\xi})_T$, where
$(U_{\xi})_T=U_{\xi}\cap (\Om G)_T$ and
$(\Om_{\xi})_T=\Om_{\xi}\cap (\Om G)_T=\Om_{\xi}$.
{}From Proposition 2.4 we deduce:

\proclaim{Proposition 5.3}
The fibre of 
$(U_{\xi})_T\to\Om_{\xi}$ over $\ga_{\xi}$ is 
$$
(\La_{e,\alg}^+\Gc)_T\cdot\ga_{\xi}=
\exp\,(\u_{\xi})_T\cdot \ga_{\xi} \cong \exp\,(\u_{\xi})_T 
\cong (\u_{\xi})_T
$$
where $(\u_{\xi})_T = 
{\tsize\bigoplus}_{0<2i<r(\xi)}\la^{2i}(\f^{\xi}_{2i})^\perp$.
\qed
\endproclaim

\no As in Proposition 2.7, we have a similar description of a 
Zariski open subspace
$$
(\La_{N',\alg}^+\Gc)_T\cdot\ga_{\xi}=
\exp\,(\u^0_{\xi})_T\cdot \ga_{\xi}\cong \exp\,(\u^0_{\xi})_T
\cong (\u^0_{\xi})_T
$$
of $(U_{\xi})_T$, where $(\u^0_{\xi})_T = 
{\tsize\bigoplus}_{0\le 2i<r(\xi)}\la^{2i}(\f^{\xi}_{2i})^\perp$.

We may now proceed as in \S 4. If $\Phi:M\to (\Om_{\alg}^k G)_T$ 
is an extended solution, then there exists some $\xi\in I'$, 
and some discrete subset $D$ of $M$, such that 
$\Phi(M-D)\sub (U_{\xi})_T$. Moreover, by taking a larger $D$ if
necessary, we may assume that 
$\Phi(M-D)\sub \exp\,(\u^0_{\xi})_T\cdot \ga_{\xi}$, and so
$$
\Phi\vert_{M-D}=A\cdot \ga_{\xi},\quad A=\exp\,C
$$
where $C:M\to (\u^0_{\xi})_T$ is meromorphic. If $\phi:M\to\sqrt e$
is the corresponding harmonic map, i.e. $\phi(z)=\Phi(z,-1)$,
then the connected component of $\sqrt e$ containing the image
of $\phi$ must be the symmetric space $N_{\xi}$. 

In Theorem 4.5 we demonstrated the existence of an element 
$\ga\in\Om G$ such that $\ga\Phi(M-D)\sub U_{\xi_{\can}}$, where
$\xi=\sum_{t} n_{t}\xi_{t}$ and $\xi_{\can}=
\sum_{t}\xi_{t}$. (As in \S 4, we write $\xi_1,\dots,\xi_l$
for the elements of $I'$ which are dual to the simple roots
$\al_1,\dots,\al_l$.) In the presence of the additional
condition $T(\Phi)=\Phi$, we can make a stronger statement:

\proclaim{Theorem 5.4} 
Assume that $G$ is semisimple, with trivial centre.
Let $\Phi:M\to (\Om_{\alg}^k G)_T$ be an extended solution 
(of finite uniton number). As explained above, we have 
$\Phi(M-D)\sub (U_{\xi})_T$ for some $\xi=\sum_{t} n_{t}\xi_{t}$.
Then there exists some $\ga \in \Om_{\alg}^k G$ with $\ga(-1)=e$,
such that $\ga\Phi(M-D)\sub (U_{\xi'_{\can}})_T$, where
$\xi'_{\can}=\sum_{n_{t}\ \text{odd}}\xi_{t}$.
\endproclaim

\demo{Proof}
Let $\hat\xi=\xi-\xi'_{\can}$.  We shall prove that
$A\cdot \ga_{\hat\xi}$ is both holomorphic and anti-holomorphic,
hence independent of $z$. The statement concerning $\ga\Phi$ will
then follow, exactly as in Theorem 4.5.
Since $A$ contains only even powers 
of $\la$, the same is true of $A^{-1}A_z=\sum_{i\ge 0}\la^i A'_i$,
i.e. we have $A'_i=0$ for all odd $i$.  As in the proof of
Theorem 4.5, we have to show that $\Im\, A'_i\sub \f^{\hat\xi}_i$ 
for all $i\ge 0$. The extended solution equation (Proposition 4.4) 
gives $\Im\,A'_i\sub \f^{\xi}_{i+1}$ for all $i\ge 0$. Thus,
it suffices to show that $\f^{\xi}_{2j+1}\sub \f^{\hat\xi}_{2j}$
for all $j\ge 0$.  We have $\f^{\xi}_{2j+1}\sub \f^{\hat\xi}_{2j+1}$
as usual. Since 
$\f^{\hat\xi}_{2j+1}=\f^{\hat\xi}_{2j}\oplus \g^{\hat\xi}_{2j+1}$,
we need to prove that $\g^{\hat\xi}_{2j+1}=0$. But this is clear,
as $\al(\hat\xi)\in 2\sq\Z$ for any positive root $\al$,
because each coefficient of $\xi_i$ in $\hat\xi$ is (by
construction) even. The required $\ga$ is given by 
$\ga^{-1}=A\cdot\ga_{\hat\xi}$.  Finally, to prove that $\ga(-1)=e$,
we observe that the (constant) extended solution $\ga^{-1}:M\to 
(U_{\hat\xi})_T$ gives a (constant) harmonic map 
$\ga^{-1}(-1):M\to N_{\hat\xi}\sub G$, and
$N_{\hat\xi}=\{ g(\exp\,\pi\hat\xi)g^{-1}\st g\in G\}=\{e\}$.
\qed\enddemo

\no It should be noted that $N_{\xi}=N_{\xi'_{\can}}$ here.
In view of Proposition 5.2 (2), 
we obtain the following upper bound on the minimal
uniton number:

\proclaim{Corollary 5.5}
Assume that $G$ is semisimple, with trivial centre.
Let $N$ be a (compact) inner symmetric space, embedded in $G$ as
a connected component of $\sqrt e$.
Let $\phi:M\to N$ be a harmonic map of finite uniton 
number. Then the minimal uniton number of $\phi$
is not greater than
$r(N)=\max\{r(\xi)\st \xi\ \text{canonical}, N_{\xi}=N\}$.
\qed
\endproclaim

In the terminology of \cite{Bu-Ra}, $r(N)$ is the maximum
height of any flag manifold which fibres canonically over $N$.
The computation of $r(N)$ is a purely Lie-algebraic matter
(see part C of Chapter 4 of \cite{Bu-Ra}).  The following 
list contains $r(N)$ for all compact irreducible inner 
symmetric spaces of classical type.

\bigpagebreak

\settabs\+\indent  &  NNNNNNNNNNNNNNNNNNNNNNNNNNN\quad &  &\cr

\+  & $N$ & $r(N)$ &\cr

\+ & -------- & -------- &\cr

\+ & ${SU_n}/{S(U_m\times U_{n-m})},\ 1\le m<n/2$ & $2m$ &\cr

\+ & ${SU_n}/{S(U_m\times U_{n-m})},\ 1\le m=n/2$ & $2m-1$ &\cr

\+ & ${SO_{2n+1}}/{SO_{m}\times SO_{2n+1-m}},\ 1\le m<n$ & $2m$ &\cr

\+ & ${SO_{2n+1}}/{SO_{m}\times SO_{2n+1-m}},\ 1\le m=n$ & $2m-1$ &\cr

\+ & ${Sp_n}/{Sp_m\times Sp_{n-m}},\ 1\le m<n/2$ & $4m$ &\cr

\+ & ${Sp_n}/{Sp_m\times Sp_{n-m}},\ 1\le m=n/2$ & $4m-2$ &\cr

\+ & ${Sp_n}/{U_n},\ n\ge 1$ & $2n-1$ &\cr

\+ & ${SO_{2n}}/{SO_{2m}\times SO_{2n-2m}},\ 2\le 2m\le n-2$ & $4m$ &\cr

\+ & ${SO_{2n}}/{SO_{2m}\times SO_{2n-2m}},\ 2\le 2m=n-1$ & $4m-1$ &\cr

\+ & ${SO_{2n}}/{SO_{2m}\times SO_{2n-2m}},\ 2\le 2m=n$ & $4m-3$ &\cr

\+ & ${SO_{2n}}/{U_n},\ n\ge 3$ & $2n-4$ &\cr

Notable examples are the sphere $S^{2n}$ ($r=2$ if $n\ge 2$), 
complex projective space $\C P^n$  ($r=2$ if $n\ge 2$), 
quaternionic projective space $\H P^n$ ($r=4$ if $n\ge 2$),
and the complex quadric $Q_n=SO_{n+2}/SO_2\times SO_n$
($r=4$ if $n\ge 5$). The result for the Grassmannian
$Gr_m(\C^n)={SU_n}/{S(U_m\times U_{n-m})}$ was conjectured
by Uhlenbeck (Problem $9$ of \cite{Uh}), and independent
proofs have been given recently in Chapter 20 of \cite{Gu}
and in \cite{Do-Sh}.

In \S 4, we saw that harmonic maps into Lie groups with
$r(\phi)\le2$ are particularly easy to describe.  For maps 
into symmetric spaces, we have similar results for $r(\phi)\le 3$.
We begin with $r(\phi)=2$.

\proclaim{Corollary 5.6}
If $\Phi:M\to (\Om^k_{\alg}G)_T$ is an extended solution
with $r(\Phi)\le 2$, then
$\Phi$ is an $S^1$-invariant extended solution.
\endproclaim

\demo{Proof}
By Proposition 5.3, we have $(U_{\xi})_T=\Om_{\xi}$ if $r(\xi)\le2$.
\qed
\enddemo

\no In particular, if $N$ is an inner symmetric space with $r(N)=2$,
then every harmonic map $\phi:S^2\to N$ is given by the twistor
construction. As we have noted, this is the case for $N=S^{2n}$
or $\C P^n$. Therefore, Corollary 5.7 gives the well known
theorem of Calabi (for $S^{2n}$) and Eells and Wood (for $\C P^n$).

For $r(\phi)=3$, we have
the following analogue of Corollary 4.9:

\proclaim{Corollary 5.7}
Let $\xi\in I'$ with $r(\xi)\le 3$.
A holomorphic map $\Phi:M-D\to (U_{\xi})_T$ is an extended solution
if and only if $u_{\xi}\circ\Phi:M-D\to \Om_{\xi}$ is
super-horizontal. 
\endproclaim

\demo{Proof} The conditions for $\Phi$ to be an extended
solution are $\Im A'_0\sub \f^{\xi}_1$, $\Im A'_1\sub \f^{\xi}_2$. 
The first of these is the condition 
for $\Phi$ to be super-horizontal; the second is vacuous as $A'_i=0$
for all odd $i$.
\qed\enddemo

\no As in Corollary 4.9, we conclude that harmonic maps
into symmetric spaces with $r(\phi)=3$ correspond to
pairs $(\Phi,\sigma)$, where $\Phi$ is a holomorphic 
super-horizontal map into $\Om_{\xi}$ and $\sigma$ is a
meromorphic section of the holomorphic vector bundle
$\Phi^\ast (U_{\xi})_T$.

With obvious modifications, the remarks in \S 4 on Weierstrass
representations, transforms, and factorizations apply equally
to harmonic maps into symmetric spaces.

\newpage
\head
Appendix A: Harmonic maps of low uniton number
\endhead

We summarize here the various special cases of our
results on harmonic maps from $S^2$ to a Lie group $G$
or symmetric space $N$, which we encountered in \S 4 and \S 5.
There are three essentially different kinds of behaviour.

{\it \no a) Twistor construction}

The twistor construction of \cite{Bu-Gu} gives all harmonic
maps which arise from $S^1$-invariant extended solutions;
the latter are simply holomorphic super-horizontal maps
into generalized flag manifolds.
This includes maps with $r(\phi)=0$ or $1$ (constant maps or
holomorphic maps into Hermitian symmetric spaces, respectively).
It also includes all harmonic maps $S^2\to N$ with $r(\phi)=2$
(Corollary 5.6).
For $N=S^{2n}$ or $\C P^n$, this means all harmonic maps.

{\it \no b) Twistor construction $+$ meromorphic section}

Any harmonic map $S^2\to G$ with $r(\phi)\le 2$, or
any harmonic map $S^2\to N$ with $r(\phi)\le 3$, is given
by a pair $(\theta,\sigma)$, where $\theta$ is a harmonic map
obtained via the twistor construction, and $\sigma$ is a
meromorphic section of a holomorphic vector bundle (over
a generalized flag manifold).  For precise statements, see
Corollaries 4.9 and 5.7. For $G=SU_3$, $N=Sp_2/U_2\cong Q_3$
or $Gr_2(\C^4)\cong Q_4$, all harmonic maps arise this way.

{\it\no c) Twistor construction $+$ meromorphic section(s)
satisfying first order differential equations}

This is the general case.  The differential equations
are solvable, locally, by integration.  The example
$G=SU_4$, which we gave in \S 4, shows how simple
these equations can be, however.  This is typical of the
situation for harmonic maps $S^2\to G$ with $r(\phi)=3$,
or harmonic maps $S^2\to N$ with $r(\phi)=4$ --- i.e. the
first level beyond b).  All harmonic
maps are of this form for the following families of 
symmetric spaces: $Gr_2(\C^n)$ for $n\ge 5$; $Q_n$ for $n\ge 5$;
$\H P^n$ for $n\ge 2$.  The same is true for the isolated
groups $G=SU_4$, $SO_5$, $Sp_2$, or $SO_6$ (and
the symmetric space $N=SO_{8}/U_4\cong Q_6$).

\newpage
\head
Appendix B: The Birkhoff decomposition
\endhead

So far, we have not made any use of the Birkhoff decomposition
of $\Om G$.  For extended solutions of finite uniton number,
it is the (finite dimensional) Bruhat manifolds which seem
to be most relevant. On the other hand, a useful consequence
of the Birkhoff decomposition is the existence of a \ll big
cell\rrr, namely the Birkhoff manifold $\La^-\Gc\cdot e$,
which is an open dense subspace of (the identity component
of) $\Om G$. This gives a natural coordinate chart on $\Om G$.
For example, in \cite{Do-Pe-Wu}, a special
role is played by extended solutions $\Phi:M\to\Om G$ whose
image lies in the big cell.  In addition to the extended solution
condition $\Im\,\la\Phi^{-1}\Phi_z\sub\La^+\gc$, we also
have in this case the condition 
$\Im\,\Phi^{-1}\Phi_z\sub\La^-\gc$, so we conclude that
$\Phi^{-1}\Phi_z=\tfrac1\la V$ for some meromorphic function 
$V:M\to\gc$. (More precisely, we conclude that 
$\Phi^{-1}\Phi_z$ is of the form $U+\tfrac1\la V$; but
$U$ is necessarily zero because $\Phi(z,1)=e$ for all $z$.)
This $V$ is called the \ll Weierstrass data\rr in \cite{Do-Pe-Wu}.

After translation by a suitable element of $\Om G$,
any extended solution has this special form, away from a
discrete subset.  In the case of the Weierstrass representation
of \S 4 and \S 5, this may be accomplished explicitly
because of the following fact:

\proclaim{Proposition B1}  Let $U_{\xi}=
\La^+_{\alg}\Gc\cdot \ga_{\xi}$ (as in \S 4). Then a dense
open subset of $\ga_{\xi}^{-1}U_{\xi}$ is contained in 
the big cell of $\Om G$.
\endproclaim

\demo{Proof} Recall (from Proposition 2.7) that we have a
dense open subset $\exp\,\u^0_{\xi}\cdot\ga_{\xi}$ of $U_{\xi}$.
An element $C$ of $\u^0_{\xi}$ is of the form
$C=C_0+\la C_1+\dots+\la^{r-1}C_{r-1}$, where
$C_i=c^{i+1}_i+ c^{i+2}_i+\dots+c^r_i$, and 
$c^j_i\in\g^{\xi}_j$. We have $\Ad\ga_{\xi}^{-1} c^j_i=
\la^{-j}c^j_i$.  Hence 
$\ga_{\xi}^{-1}(\exp\,\u^0_{\xi}\cdot\ga_{\xi})=$ {}
$\exp\,\Ad\ga_{\xi}^{-1}\u^0_{\xi}\cdot e
\sub$ {} $\La^-\Gc\cdot e$, as required
\qed
\enddemo

For example, consider the extended solution $\Phi=\exp\,C\cdot\ga$
(for a harmonic map $S^2\to U_4$) from \S 4, where
$$
C=
\pmatrix
{0} & a_1 & b_1 & c_1\\
{0} & {0} & a_2 & b_2\\
{0} & {0} & {0} & a_3\\
{0} & {0} & {0} & {0}
\endpmatrix
+\la 
\pmatrix
{0} & {0} & d_1 & e_1\\
{0} & {0} & {0} & d_2\\
{0} & {0} & {0} & {0}\\
{0} & {0} & {0} & {0}
\endpmatrix
+\la^2
\pmatrix
{0} & {0} & {0} & f_1\\
{0} & {0} & {0} & {0}\\
{0} & {0} & {0} & {0}\\
{0} & {0} & {0} & {0}
\endpmatrix,\ 
\ga=
\pmatrix
\la^3 & {} & {} & {}\\
{} & \la^2 & {} & {}\\
{} & {} & \la & {}\\
{} & {} & {} & 1
\endpmatrix.
$$
A calculation shows that
$$
\Phi^{-1}\Phi_z=
\tfrac1\la
\pmatrix
{0} & a'_1 & d'_1 & f'_1\\
{0} & {0} & a'_2 & d'_2\\
{0} & {0} & {0} & a'_3\\
{0} & {0} & {0} & {0}
\endpmatrix
$$
in this case.
Thus, the Weierstrass data consists precisely of the 
derivatives of the rational functions 
$a_1,a_2,a_3,d_1,d_2,f_1$. We saw in \S 4 that these 
rational functions (locally) parametrize extended 
solutions of this type.

$${}$$

\eightpoint
 
\no School of Mathematical Sciences, University of Bath, 
Bath BA2 7AY, United Kingdom\newline
Department of Mathematics, University of
Rochester, Rochester, New York 14627, USA

\no f.e.burstall\@maths.bath.ac.uk\newline
gues\@db1.cc.rochester.edu

\enddocument